\documentclass[11pt]{article}

\usepackage[T1]{fontenc}
\usepackage{lmodern} % (optional) good default fonts
\usepackage{newtxtext}
\usepackage[varvw]{newtxmath}

\usepackage[margin=1in]{geometry}
\usepackage{xcolor}
\usepackage{graphicx}
\usepackage[bottom]{footmisc}
\usepackage{multicol} % keep if you use it later
\usepackage{makeidx}  % keep if you use an index later
\makeindex

\usepackage{mathrsfs}
\usepackage{bm}
\usepackage{amsmath}

\usepackage{algorithm}
\usepackage{algpseudocode}

\usepackage{soul}
\usepackage{cancel}

\usepackage[numbers]{natbib}

\usepackage[hidelinks]{hyperref}

\def\Cc{{\mathscr C}}
\def\NN{{\mathbb{N}}}

\def\CC{{\mathbb{C}}}
\def\RR{{\mathbb{R}}}

\def\Qq{{\mathscr Q}}
\def\Pp{{\mathscr P}}
\def\E{{\mathcal E}}
\def\L{{\mathcal L}}

\def\func#1#2{ \colon {#1} \longrightarrow {#2}}

\def\epsilon{\varepsilon}

\newcommand{\mw}[1]{\text{\sffamily\normalshape #1\/}}
\newcommand{\eff}{\mw{eff}\kern1pt}
\newcommand{\tol}{\mw{tol}\kern1pt}
\newcommand{\Sparsify}{\mw{Sparsify}\kern1pt}
\newcommand{\rms}{\mw{rms}\kern1pt}
\def\network\_recovery{\textsf{network\_recovery}\kern1pt}

\usepackage{amsthm}
\newtheorem{propo}{Proposition}[section]

\newtheorem{example}[propo]{Example}
\newtheorem{remark}[propo]{Remark}
\newtheorem{problem}[propo]{Problem}
\newtheorem{definition}[propo]{Definition}

\newtheorem{theorem}[propo]{Theorem}

\title{Simultaneous recovery of a sparse topology\\ and the admittance of an electrical network}

\author{\'Alvaro Samperio}
\date{} % empty date (like many Springer chapters)

\begin{document}

\maketitle

% \begingroup
% \renewcommand\thefootnote{}
% \footnotetext{Departamento de Matem\'aticas, CUNEF Universidad, Madrid, Spain.\\
%  Email: \texttt{alvaro.samperio\symbol{'100}cunef.edu}\\The author was funded by a FPI grant of the Research Project PGC2018-
%  096446-B-C21 (with the help of the FEDER Program).
%  The author was partially supported by the project PID2022-138906NB-C21 funded by
%  MICIU/AEI/ 10.13039/501100011033 and by ERDF/EU. The author was partially supported by
%  the project 2024PROD00092 funded by AGAUR. \\
%  Math. Subject Class. 2020: 94C15, 05C85, 05C90. Keywords: graph sparsification, ill-posed inverse problems, electrical networks, sparse linear regression.
%  }
% \endgroup

\begingroup
\renewcommand\thefootnote{}
\footnotetext{Departamento de Matem\'aticas, CUNEF Universidad, Madrid, Spain.\\
 Email: \href{mailto:alvaro.samperio@cunef.edu}{\texttt{alvaro.samperio\symbol{'100}cunef.edu}}
 }
\endgroup

\begin{abstract}
We show that the problem of recovering the topology and admittance of an electrical network from power and voltage data at all vertices is often ill-posed, and sometimes it even has multiple solutions. We reformulate the problem to seek for a sparse network, {\it i.e.}, with few edges, which fits the data up to a given tolerance. We propose an algorithm to solve this reformulated problem. It combines, in an iterative procedure, the resolution of non-negative linear regression problems, and techniques of spectral graph sparsification. The algorithm is based on original results bounding the fitting error of a sparse approximation of a network. We illustrate our techniques with several experimental results in which we are able to recover a sparse network.
\end{abstract}

%%%%%%%%%%%%%%%%%%%%%%%%%%%%%%%%%%%%%%%%%%%%%%%%%%%%%%%%%%%%%%%%%%%%%%%%%%%%%%%%%%%%%%%%% c3

\section{Introduction and preliminaries}\label{sec:intro}

We begin with some notations, definitions and results of vector calculus on electrical networks from \cite{S24}. Given $m\in \NN^*=\NN\setminus \{0\}$ and a finite set $V$, we denote the spaces of real and complex vector functions on $V$ as $\Cc(V,\RR^m)=\{u\func{V}{\RR^m}\}$ and $\Cc(V,\CC^m)=\{u\func{V}{\CC^m}\}$, respectively. In the particular case of real scalar functions, we denote $\Cc(V)=\Cc(V,\RR)$. We denote the set of nonnegative scalar functions on $V$ as $\Cc^+(V)$. Given $\bm u\in \Cc(V,\CC^m)$, respectively $\bm u \in \Cc(V,\RR^m)$, the components of $\bm u$ are $u_1,\ldots,u_m\in \Cc(V,\CC)$, respectively $u_1,\ldots,u_m\in \Cc(V)$, {\it i.e.} $\bm u=(u_1,\ldots,u_m)$. We denote its real and imaginary parts as $\Re(\bm u)=(\Re(u_1),...,\Re(u_m))\in \Cc(V,\RR^m)$ and $\Im(\bm u)=(\Im(u_1),...,\Im(u_m))\in \Cc(V,\RR^m)$, and its complex conjugate as $\overline{\bm u}=(\overline{u_1},\ldots,\overline{u_m})$.

 Given $F\subseteq V$,  we denote by $\Cc(F)$ the subspace of the functions in $\Cc(V)$ which are null on $V\setminus F=F^c$. Also, for each $u\in \Cc(V,\CC)$ we define
 $\int_Fu dx=\sum\limits_{x\in F}u(x)$. On $\Cc(V,\CC^m)$ we consider the inner product 
$$\langle \bm u,\bm v\rangle=\sum_{j=1}^{m}\langle u_j,v_j\rangle=\sum_{j=1}^{m}{\int_Vu_j\bar v_jdx},$$
whose associated norm is denoted by $||\cdot||_2$, or just by $||\cdot||$. This inner product induces an inner product on $\Cc(V,\RR^m)$. We also consider  in $\Cc(V,\CC^m)$ the norm $||\cdot||_{_\infty}$, which is defined for each $\bm u=(u_1,...,u_m)\in\Cc(V,\CC^m)$ as $||\bm u||_{_\infty}=\max\limits_{j=1,..,m}\left\{||u_j||_{_\infty}\right\}$.

A Direct Current (DC) network is modeled mathematically as a graph with a real positive weight at each edge, which is the inverse of its {\it resistance}. For modeling Alternating Current (AC) networks there are different alternatives, depending on the characteristics of the network. In this work, we consider balanced three-phase networks in which all lines are inductive and ``short'', ({\it i.e.}, their length is shorter than 80km). In that case, the {\it susceptance to earth} of the lines can be neglected and the mathematical model of a network is a graph in which each physical connection between a pair of vertices is represented by a single edge characterized by a complex weight with nonnegative real part and nonpositive imaginary part called {\it admittance}, (see \cite{SE68, W12}). 

\begin{definition}\label{defred}{\rm
	An {\it electrical network} is a pair $\Gamma = (V, a)$ where $V$ is a finite nonempty set called {\it vertex set}, and $a\in \Cc(V\times V, \CC)$, called {\it admittance}, satisfies the following properties:
 	\begin{enumerate}%[{\rm (i)}]
	\item $a=c-i b$ where $c,b\in \Cc^+(V\times V)$ and the function $c$ is called {\it conductance} whereas the function $b$ is called {\it susceptance}.
	\item $c$ and $b$ are symmetric and $c(x,x)=b(x,x)=0$, for any $x\in V$. 
	\end{enumerate}}
\end{definition}

An electrical network $\Gamma = (V, a)$ has an associated graph $G(\Gamma)=(V, E(\Gamma))$ called its ``{\it network topology}''. Its vertex set is $V$ and its {\it edge set} is $E(\Gamma) =\big\{\{x,y\} \, \, \text{such that}  \, \,  x,y\in V \, \, \text{and}  \, \,  a(x,y)\neq 0 \big\}$. We say that $x,y\in V$ are {\it adjacent} iff $\{x,y\}\in E(\Gamma)$, and we denote $e_{xy}=e_{yx}=\{x,y\}$. For any $e_{xy}\in E(\Gamma)$, we denote by $a(e_{xy})$ the value $a(x,y)=a(y,x)\neq 0$. Note that $a(x,x)=0$ for any $x\in V$, so $(V, E(\Gamma))$ is a finite undirected simple graph without loops, {\it i.e.}, without edges from a vertex to itself.

In the particular case the network is operating in DC, $b = 0$, so the network $\Gamma = (V, c)$ is a real weighted graph.

For an AC network $\Gamma = (V, a)$ we consider the two DC networks $\Gamma_c = (V, c)$ and $\Gamma_b = (V, b)$, which we will call the conductance and susceptance networks of $\Gamma$, respectively. Note that $E(\Gamma)=E(\Gamma_c) \cup E(\Gamma_b)$.

For any function $a\in \Cc(V\times V, \CC)$  that is symmetric and satisfies that $a(x,x)=0$ for any $x\in V$, we denote by $\L_a$ the endomorphism of $\Cc(V,\CC)$ defined as 
$$\L_a(u)(x)\displaystyle=\sum\limits_{y\in V}
a(x,y)\,\big(u(x)-u(y)\big)
$$
for any $u\in \Cc(V, \CC)$ and any $x\in V$. When $a$ is the admittance of an electrical network $\Gamma=(V,a)$, $\L_a$ is called the {\it Laplacian} of $\Gamma$. {\it Gauss' Theorem} states that for any $u\in \Cc(V, \CC)$, $\int_{V} \L_a(u) dx=0$.

If $\Gamma=(V,a)$ is an AC network, at a given time, the potential at the vertices can be represented by $u\in \Cc(V, \CC)$, the current injected at the vertices is equal to $\L_a(u)\in \Cc(V, \CC)$ and the power injected at the vertices is equal to $s=u \overline{\L_a(u)}\in \Cc(V, \CC)$. Similarly, if $\Gamma=(V,c)$ is a DC network, the potential at the vertices can be represented by $u\in \Cc(V)$, the current injected at the vertices is equal to $\L_c(u)\in \Cc(V)$, and the power injected at the vertices is equal to $s=u \L_c(u)\in \Cc(V)$. These equations relating voltage and power are called the {\it power flow equations}, (see \cite{montes}).

\begin{remark}\label{vmin}
	{\rm In order to avoid possible faults in the operation of an electrical network,
		voltages are usually maintained close to a reference value, that is chosen to be
		$u=1$ in appropriate units (see \cite{montes}). In this chapter we will fix minimum and
		maximum values of voltage $0<|u|_{\rm min}\le |u|_{\rm max}$. Later, in the experiments, we will use $|u|_{\rm min}=0.9$ and $|u|_{\rm max}=1.1$.}
\end{remark}

Let $\Gamma=(V,a)$ be an electrical network and let $F\subsetneq V$. We say that a function $u$ is harmonic on $F$ when $\L(u)=0$ on $F$. Given $g\in \Cc(F^c,\CC)$, the {\it Dirichlet problem} of finding $u\in \Cc(V,\CC)$ such that
	\begin{equation}\label{eqprboundary0} \L_a(u)=0\hspace{.25cm}\hbox{on}\hspace{.25cm}
		F,\hspace{.5cm}
		u=g\hspace{.25cm}\hbox{on}\hspace{.25cm} F^c,
	\end{equation}
always has a solution. Moreover, the images of all the solutions to (\ref{eqprboundary0}) by the Laplacian $\L_a$ are equal. Furthermore, there is an endomorphism $\Lambda$ of $\Cc(F^c,\CC)$ called the {\it Dirichlet-to-Neumann map of $\Gamma$ and $F$} defined for each $g\in \Cc(F^c,\CC)$ as $\Lambda(g)=\L_a(u_g)$, where $u_g$ is a solution to (\ref{eqprboundary0}), (see \cite{S24}). 

That is, the Dirichlet-to-Neumann map is the linear map which gives us the relationship between the potential at $F^c$ and its corresponding current injected at $F^c$ under the condition that the potential is harmonic on $F$.

Additionally, for all DC networks and for some AC networks, the Dirichlet-to-Neumann map of $\Gamma$ and $F$ is the Laplacian of another network whose vertex set is $F^c$ and whose admittance we denote by $a^{\Lambda}$, which is called the {\it Kron reduction of $\Gamma$ with respect to $F$}. When $x$ and $y$ are two distinct vertices of $V$, the Dirichlet-to-Neumann map of $\Gamma$ and $F=V \setminus\left\{x,y\right\}$ is always a Kron reduction (with two vertices), and the value $a^e(x,y):=a^{\Lambda}(x,y)$ is called the {\it effective admittance between $x$ and $y$}. When $\Gamma$ is a DC network, the effective admittance is nonnegative, it is called effective conductance and it is denoted as $c^e(x,y)$.

Denoting by $c_{\setminus e_{xy}}^{e}(x,y)$ the effective conductance between $x$ and $y$ in the network obtained from a DC network $\Gamma$ by removing the edge $e_{xy}$, we have that, (see \cite{S24}),
\begin{equation}\label{eqsumadc}
c^{e}(x,y)=c_{\setminus e_{xy}}^{e}(x,y)+ c(x,y),
\end{equation}
and thus $c^{e}(x,y)\geq c(x,y)$. For every edge $e_{xy}$ we define the {\it effective resistance between $x$ and $y$} as $r^{e}(x,y)=1/c^{e}(x,y)>0$.

We define the {\it energy} of a DC network $\Gamma=(V,c)$, (see \cite{C18}), as the positive semidefinite symmetric bilinear form on $\Cc(V)$ given by 
$$\displaystyle \E(u,v)=\langle u, \L_c(v) \rangle=\int_Vu \L_c(v)dx.$$

Nowadays, electrical systems are expanding and increasing in complexity really fast due to the deregulation and proliferation of distributed energy resources. Because of this expansion, sometimes the system operators do not have precise and updated
information about the admittance or the topology of the network. Also, the dependence of admittance on temperature can cause the information about it to be incorrect \cite{rev}. The lack of topology information is especially frequent in the secondary distribution network \cite{quitabajac}, in which the lines are short \cite{distribucion_corta}. Moreover, recently, the availability and accuracy of measured data in the electrical
system has increased significantly. At the nodes, power injection and voltage magnitude can be measured using smart meters and the phase angle information can be obtained by micro phasor measurement units \cite{quitabajac}. In that scenario, the simultaneous recovery of the topology and admittance of a network from voltage and power data is of great interest from an applied point of view.

Several authors have formulated different inverse problems on networks which involve the recovery of the topology, starting from different data; and they have proposed methods to solve these problems. For instance, in the context of the inverse conductance problem, in \cite{GeRo22} the authors introduced a method which recovers the conductance and the topology up to vertices of combinatorial degree two of any DC tree from its Dirichlet-to-Neumann map.

In \cite{rev}, the application of a least squares approach to solve two problems is discussed: the first one concerns the admittance estimation alongside voltage estimation from a known topology, and the second one consists of the voltage estimation alongside detection of errors in the network topology determined by a given admittance. However, in that paper, the simultaneous recovery of topology and admittance is not discussed. 

An iterative method is proposed in \cite{quitabajac} to recover simultaneously the topology and admittance of a network from power and voltage data, starting at each step from an estimated admittance and then calculating the current flow in each edge, removing an edge if its conductance and susceptance have a value below a given threshold. Nevertheless, relying only on the values of the admittance to decide if an edge can be removed from a network can lead us to an error, as we show in Example \ref{ej3}. The method proposed in \cite{aproximacion_dc} also estimates both the topology and admittance of a network from power and voltage data, but only in the particular case of distribution networks which operationally are a tree.

In this work, we study the problem of simultaneously recovering the topology and admittance of an
electrical network from power and voltage data at all vertices. We introduce the following notation to formulate the general network recovery problem with exact data (Problem \ref{fijoexacto}).

For any $a\in \Cc(V\times V, \CC)$  that is symmetric and satisfies that $a(x,x)=0$ for any $x\in V$, we denote by $\Pp_{a}\func{\Cc(V,\CC^m)\times\Cc(V,\CC^m)}{\Cc(V,\CC^m)}$ the function defined for each $\left(\bm v,\bm w\right)=\left((v_1,...,v_m),(w_1,...,w_m)\right)\\ \in \Cc(V,\CC^m)\times\Cc(V,\CC^m)$ as 
$$\Pp_{a}\left(\bm v,\bm w\right)=(v_1\overline{\L_a(w_1)},...,v_m\overline{\L_a(w_m)}).$$
Of course, when $a$ is a real function, $\Pp_{a}\left(\bm v,\bm w\right)\in \Cc(V,\RR^m)$ for every $\bm u, \bm v\in \Cc(V,\RR^m)$. Also, if $a$ is the admittance of an AC network, respectively a DC network, then for any $\bm u\in \Cc(V,\CC^m)$, respectively $\bm u\in \Cc(V,\RR^m)$, the j-th component of the function $\Pp_{a}\left(\bm u,\bm u\right)$ is equal to the power injected at the network when the potential is $u_j$, for any $j=1,...,m$.

We denote as $E(K_{V})$ the set of edges of the complete graph on $V$, {\it i.e.} the set of all possible edges on the vertex set $V$. In an applied situation, one of the main causes that the system operators of a network $\Gamma$ may not know its topology $(V,E(\Gamma))$ is the existence of lines which have a switch whose status is unknown, so they do not know whether their corresponding edges belong to $E(\Gamma)$ or not, (see \cite{rev}). In general, we formulate the inverse problem of recovering the topology and admittance
of $\Gamma$ with the {\it a priori} information that there is a set $E\subseteq E(K_{V})$ of edges which are candidates to belong to $E(\Gamma)$, and we know that the edges in $E(K_{V})\setminus E$ do not exist in the network topology. Note that, the particular case in which $E=E(K_{V})$, corresponds with the framework in which we do not have any {\it a priori} information about $E(\Gamma)$.

\begin{problem}\label{fijoexacto}
	Let $V$ be a finite set of vertices, let $E$ be a set of edges on $V$, let $0<|u|_{\rm min}\le |u|_{\rm max}$ be positive values, let $m\in \NN^*$ and let $\bm u,\bm s\in \Cc(V,\CC^m)$, respectively $\bm u,\bm s\in  \Cc(V,\RR^m)$, such that $|u|_{\rm min}\leq |u_j| \leq |u|_{\rm max}$ for all $j=1,\ldots,m$. Determine an AC, respectively a DC, electrical network $\Gamma=(V,a)$ with set of edges $E(\Gamma)\subseteq E$ such that $\bm s=\Pp_{a}(\bm u,\bm u)$.
\end{problem}

The pair $(\bm u,\bm s)$ is usually called {\it the data pair} or {\it the data set} and hence we always understand that for any $j=1,\ldots,m$, $u_j$ represents a measurement of voltage and $s_j$ represents its corresponding measurement of power. Therefore for any $j=1,\ldots,m$ we have that $u_j(z)\not=0$ for all $z\in V$, since $0<|u|_{\rm min}\leq |u_j| $. Note that the set of admittances, respectively conductances, on $V\times V$ that define a network $\Gamma$ satisfying $E(\Gamma)\subseteq E$ can be identified with $\Cc^+(E)\times\Cc^+(E)$, respectively $\Cc^+(E)$.

In Section \ref{s3} we show that Problem \ref{fijoexacto} is ill-posed. Therefore, we usually obtain a solution whose set of edges is the whole $E$. If $E$ is large, the recovered network is inefficient for solving usual applied problems in electrical networks which require knowing the topology and admittance. Those applications include failure identification, power flow optimization or generation scheduling \cite{aproximacion_dc}. In Section \ref{s4} we reformulate the problem to look for a sparse network whose fitting error to the data is lower than a fixed tolerance, (see Problem \ref{pr2}). Recovering a network with those characteristics would allow the efficient and accurate resolution of those applied problems.

In Section \ref{sec33} we review the notion of spectral sparsification of a network introduced in \cite{sparse_original} by Spielman and Teng. We also review and give a novel physical interpretation to Algorithm \ref{sparsify}, which was proposed in \cite{sparse_resumen} to construct a sparse spectral approximation of a network by removing edges from it. In Section \ref{sec34}, we prove original theoretical results (see Theorems \ref{cota} and \ref{cota_ac}) that give an upper bound on the increase of the fitting error of a network in that process of removing edges.

These theoretical results are the basis for the
Algorithm \ref{selfhealing} that we propose in Section \ref{section_algor} to solve Problem \ref{pr2}. The algorithm consists in an iterative procedure of recovering a network and sparsifying it,
with an automatic exploration to find the levels of closeness of sparsification that
allow us to remove edges without increasing the error in the data above the fixed tolerance.

Finally, in Section \ref{s5} we show experimental results of the application of Algorithm \ref{selfhealing}. We see that, in some cases, we recover a
network with the actual topology and low error in the admittance. We show other
cases in which we recover a Kron reduction of the network alongside isolated nodes when the power injected at them is always equal to zero in the data set. Additionally, we show that, in some cases, the
algorithm has the advantage of providing a sparse approximation of the real network, in which the computation of the mentioned applied problems is faster than in the real network.

\section{Ill-posedness of the network recovery problem}\label{s3}

Problem \ref{fijoexacto} is ill-posed in general, often having multiple solutions. For example, in the following example, we show a case in which, under the restriction that the power injected at a node is zero at a certain vertex for all the $s_j$ functions in the data pair  $(\bm u,\bm s)$, then there are multiple solutions for any choice of $(\bm u,\bm s)$. This restriction is common in networks of the real world, in which there are vertices which are not associated to any generator or consumer of power.

\begin{example}\label{triangulo}
	{\rm	
		Let $\Gamma=(V,c)$ be the DC electrical network which is the weighted path in Figure
		\ref{camino}, with vertex set $V=\left\{x,y,z\right\}$, with set of edges $E(\Gamma)=\left\{e_{xy},e_{yz}\right\}$, whose conductance values are $c(x,y)=c(y,x)>0$ and $c(y,z)=c(z,y)>0$, respectively; and with Laplacian $\L$. 
        
		\begin{figure}[h!]
			\centering
			\includegraphics[width=130mm]{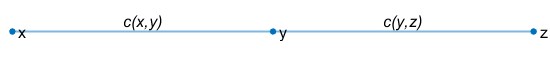}
			\caption{Electrical network in Example \ref{triangulo}.}
			\label{camino}
		\end{figure}
        
		Let $(\bm u,\bm s)$ be a data pair satisfying the power flow equations, that is, $\bm s=\Pp_{c}(\bm u,\bm u)= \\ \big(u_1\L(u_1),\ldots,u_m\L(u_m)\big)$, and moreover that $\bm s(y)=\bm 0$. We have that for each $j=1,...,m$, $u_j$ is harmonic on $y$, because $u_j(y)\neq 0$. By Gauss' Theorem, $\L(u_j)(x)=I_j$, $\L(u_j)(y)=0$ and $\L(u_j)(z)=-I_j$, for some $I_j\in\RR$. Then, each $u_j$ must be equal to $u_j(x)=\frac{I_j}{c(x,y)} +\zeta_j$, $u_j(y)=\zeta_j$ and $u_j(z)=-\frac{I_j}{c(y,z)} +\zeta_j$, for some $\zeta_j\in\RR$.
		
		We consider the Problem \ref{fijoexacto} with data pair $(\bm u,\bm s)$ and the edge set of the complete graph in $V$, $E(K_V)=\left\{e_{xy},e_{xz},e_{yz}\right\}$. Then, a DC electrical network $\Gamma'=(V,c')$ with Laplacian $\L'$ is a solution of the problem if and only if we have that $c'\in \Cc^+(E)$  is a solution of 
		$$
		\begin{cases} 
			c'(x,y)\big (u_j(x)-u_j(y)\big )+c'(x,z)\big( u_j(x)-u_j(z)\big)=I_j \\[1ex]
			c'(x,z)\big (u_j(z)-u_j(x)\big)+c'(y,z)\big (u_j(z)-u_j(y)\big)=-I_j  \\[1ex]
			c'(x,y)\big(u_j(y)-u_j(x)\big)+c'(y,z)\big(u_j(y)-u_j(z)\big)=0, 
		\end{cases}
		$$
		for all $j=1,...,m$. For each $j=1,...,m$, the third equation in that system depends linearly on the rest. If $I_j=0$, then $u_j$ is constant on $V$, so any value of the conductance $c'$ satisfies that $s_j=u_j\L'(u_j)$. If, on the contrary, $I_j\neq 0$, then the previous system is equivalent to the following system, which is independent of $I_j$ and $\zeta_j$:
		
		\begin{equation}\label{sistemaej}
			\begin{cases} 
				\dfrac{c'(x,y)}{c(x,y)}+\left(\dfrac{1}{c(x,y)}+\dfrac{1}{c(x,z)}\right)c'(x,z)=1 \\[3ex]
				\left(\dfrac{1}{c(x,y)}+\dfrac{1}{c(x,z)}\right)c'(x,z) + \dfrac{c'(y,z)}{c(x,z)}=1.
			\end{cases}
		\end{equation}
		
		The function $\alpha\in \Cc(E)$; which is determined by the values $\alpha(x,y)=-1-\frac{c(x,y)}{c(y,z)}$, ${\alpha(x,z)=1}$ and $\alpha(y,z)=-1-\frac{c(y,z)}{c(x,y)}$; is a solution to the homogeneous system associated to (\ref{sistemaej}). Therefore, any network $\Gamma'=(V,c')$ whose conductance is of the form $c'=c+\delta \alpha$, with $\delta \in \left[0, \frac{c(x,y)c(y,z)}{c(x,y)+c(y,z)}\right]$, is a solution to Problem \ref{fijoexacto}.
		
		Among the solutions with conductance of this form $c'=c+\delta \alpha$, the one with $\delta=0$ is the real network $\Gamma$. The solution $\Gamma'=(V,c')$ with $\delta=\frac{c(x,y)c(y,z)}{c(x,y)+c(y,z)}$ satisfies $c'(x,y)=c'(y,z)=0$, so it has only one edge, $E(\Gamma')=\left\{e_{xz}\right\}$. Moreover, its conductance satisfies
		\begin{equation}\label{equivalente}
			\frac{1}{c'(x,z)}=\frac{1}{c(x,y)}+\frac{1}{c(y,z)},
		\end{equation}
		so $\Gamma'$ is the union of the isolated vertex $\left\{y\right\}$ and the Kron reduction of $\Gamma$ with respect to $\left\{y\right\}$; and $c'(x,z)$ is equal to the effective conductance in $\Gamma$ between $x$ and $z$, $c^e(x,z)$. The fact that this network is a solution of the problem is because its Laplacian is the Dirichlet-to-Neumann map of $\Gamma$ and $\left\{y\right\}$, which gives us the relationship between potential and injected current at $V\setminus\left\{y\right\}=\left\{x, z\right\}$ when the potential is harmonic on $y$. Additionally, the rest of solutions with conductance of form $c'=c+\delta \alpha$, are the infinite solutions with $\delta\in \left(0, \frac{c(x,y)c(y,z)}{c(x,y)+c(y,z)}\right)$, which have the topology of $(V,E(K_V))$, the complete graph in $V$. Note that Problem \ref{fijoexacto} can be ill-posed even in cases in which the set of edges $E$ considered is the set of edges of a solution, as happens in this example.
	}	
\end{example}

In a numerical context, we usually have errors in the data pair $(\bm u,\bm s)$, and as a consequence often there is no exact solution to Problem \ref{fijoexacto}. We reformulate the problem in Problem \ref{fijo}, in order to have a solution for any data pair. We seek a network for which the power flow equations are the closest possible to be satisfied in the least squares sense. If $\Gamma=(V,c)$ is a DC network, we define the {\it error of the network $\Gamma$ in the data pair $(\bm u,\bm s)$}, $\rms(\Gamma, \bm u,\bm s)$, as the root mean square of the set of evaluations at all vertices of the residuals $u_j\L(u_j)-s_j$ of all the power flow equations corresponding with all the pairs $\left(u_j,s_j\right)$, $j=1,\ldots,m$. That is,
$$\rms (\Gamma, \bm u,\bm s)=\frac{||\Pp_{c}(\bm u,\bm u)-\bm s||}{\sqrt{m|V|}} =\sqrt{\frac{1}{m|V|}\sum_{j=1}^{m}{\int_V\left(u_j\L(u_j)-s_j\right)^2dx}}.$$

If $\Gamma=(V,a)$ is an AC network, we define the {\it error of the network $\Gamma$ in the data pair $(\bm u,\bm s)$}, $\rms(\Gamma, \bm u,\bm s)$, as the root mean square of the set of evaluations at all vertices of the real and imaginary parts of the residuals $u_j\overline{\L(u_j)}-s_j$ of all the power flow equations corresponding with all the pairs $\left(u_j,s_j\right)$, $j=1,\ldots,m$. That is,

$$\rms (\Gamma, \bm u,\bm s)=\frac{||\Pp_{a}(\bm u,\bm u)-\bm s||}{\sqrt{2m|V|}}=\sqrt{\frac{1}{2m|V|}\sum_{j=1}^{m}{\int_V\left|u_j\overline{\L(u_j)}-s_j\right|^2dx}}.$$

\begin{problem}\label{fijo}
	Let $V$ be a finite set of vertices, let $E$ be a set of edges on $V$, let $0<|u|_{\rm min}\le |u|_{\rm max}$ be positive values, let $m\in \NN^*$ and let $\bm u,\bm s\in \Cc(V,\CC^m)$, respectively $\bm u,\bm s\in  \Cc(V,\RR^m)$, such that $|u|_{\rm min}\leq |u_j| \leq |u|_{\rm max}$ for all $j=1,\ldots,m$. Determine an AC, respectively a DC, electrical network $\Gamma=(V,a)$ with set of edges $E(\Gamma)\subseteq E$ and such that the error $\rms (\Gamma, \bm u,\bm s)$ is minimum.
\end{problem}

Finding a solution to Problem \ref{fijo} is equivalent to solving a Non-Negative Least Squares (NNLS) problem, that is, a least squares problem with the restriction that the coefficients must be nonnegative, (see \cite{SH13}). In order to see that, given $\bm u\in \Cc(V,\RR^m)$ and a set of edges $E$, we define the linear
operator $\mathcal{M}_{E,\bm u}\func{\Cc(E)}{\Cc(V,\RR^m)}$ as 
$$\mathcal{M}_{E,\bm u}(c)=\big(u_1 \L_c(u_1),\ldots,u_m \L_c(u_m)\big)=\Pp_{c}(\bm u,\bm u).$$

Then, for any DC network $\Gamma=(V,c)$ with set of edges $E(\Gamma)\subseteq E$, and any data pair $(\bm u,\bm s)\in \Cc(V,\RR^m)\times \CC(V,\RR^m)$, we have that 
\begin{equation}\label{eqrms}
\rms (\Gamma, \bm u,\bm s)=\frac{||\mathcal{M}_{E,\bm u}(c)-\bm s||}{\sqrt{m|V|}}.    
\end{equation}

Similarly, given $\bm u\in \Cc(V,\CC^m)$ and a set of edges $E$, we define the linear operator \\$\mathcal{M}_{E,\bm u}\func{\Cc(E)\times \Cc(E)}{\Cc(V,\RR^{2m})}$ as
$$\mathcal{M}_{E,\bm u}(c,b)=\Big(\Re\big (u_1\overline{\L_a(u_1)}\big ),\Im\big (u_1\overline{\L_a(u_1)}\big ),\ldots, \Re\big (u_m\overline{\L_a(u_m)}\big ),  \Im\big (u_m\overline{\L_a(u_m)}\big )\Big),$$
where $a=c-ib$. 
We also consider the operator $\mathcal{S}\func{\Cc(V,\CC^m)}{\Cc(V,\RR^{2m})}$ defined as 
$$ \mathcal{S}(\bm s)=\Big(\Re(s_1),\Im(s_1),\ldots,\Re(s_m),\Im(s_m)\Big).$$
Then, for any AC network $\Gamma=(V,a)$ with set of edges $E(\Gamma)\subseteq E$, where  $a=c-ib$, and any data pair $(\bm u,\bm s)\in \Cc(V,\CC^m)\times \CC(V,\CC^m)$, we have that 
\begin{equation}\label{eqrmsac}\notag
\rms (\Gamma, \bm u,\bm s)=\frac{||\mathcal{M}_{E,\bm u}(c,b)-\mathcal{S}(\bm s)||}{\sqrt{2m|V|}}.
\end{equation}

Both for DC and AC networks, Problem \ref{fijo} consists in seeking for a nonnegative function $c\in \Cc^+(E)$ in the case of DC networks, or two nonnegative functions $c,b \in \Cc^+(E)$ in the case of AC networks, such that the norm of the difference between its image under the real linear operator $\mathcal{M}_{E,\bm u}$ and the vector function $\bm s$, or $\mathcal{S}(\bm s)$,  is minimum, so Problem \ref{fijo} is equivalent to a NNLS problem.

Any NNLS problem is a convex quadratic optimization problem, (see \cite{SH13}), so every local minimum of it is a
global minimum. We can obtain a solution $\Gamma$ to Problem \ref{fijo} calculating a minimum of its associated NNLS problem with an interior
point method. The interior point methods are among the main numerical algorithms used to solve NNLS problems,
(see \cite{ipopt_para_nnls}). We denote the numerical solution $\Gamma$ with error $\rms  \equiv \rms (\Gamma, \bm u,\bm s)$ obtained from $E$ and $(\bm u,\bm s)$ as $[\Gamma,\rms ]=\network\_recovery (E,\bm u,\bm s)$.

A sufficient condition that implies that Problem \ref{fijo} has a unique solution is that $m \geq \frac{|E|}{|V|}$ and $\kappa(\mathcal{M}_{E,\bm u})<\infty$, {\it i.e.}, the {\it condition number of $\mathcal{M}_{E,\bm u}$} is finite, (see \cite{SH13}). We will require
in Section \ref{s5} that the data sets used to test the sparse network recovery algorithm satisfy $m \geq \frac{|E|}{|V|}$. Nevertheless, often the condition $\kappa(\mathcal{M}_{E,\bm u})<\infty$ is not satisfied, and thus there are multiple solutions to Problem \ref{fijo}; as seen in Example \ref{triangulo}, in which there are multiple solutions with zero error.

As we will see in Section \ref{s5}, in many cases in which the pair  $(\bm u,\bm s)$ contains voltage and power data corresponding to an electrical network with some error, the value of the condition number $\kappa(\mathcal{M}_{E,\bm u})$ is finite but it is very high, so Problem \ref{fijo} has a unique solution but it is severely ill-posed.

%-----------------------------------------------------------------------------------------------------------------------------------------------------------------------------------------------------------------------------------------------------------------------------------------------------------------------------

\section{Reformulation of the problem: Recovery of a sparse electrical network}\label{s4}

As we have discussed in the previous section, given a data data pair $(\bm u,\bm s)$ and a set of edges $E$, the Problem \ref{fijo} of recovering a network $\Gamma$ with minimum error such that $E(\Gamma)\subseteq E$ is ill-posed. Because of that, even if there is a solution $\Gamma^*$ to the exact Problem \ref{fijoexacto} whose topology is much more sparse than $(V,E)$, {\it i.e.}, $|E(\Gamma^*)|\ll |E|$, solving Problem \ref{fijo} we usually get a solution whose topology is $(V,E)$. If the number of edges in $E$ is high, a solution with that topology is not efficient for applications. We are
interested in recovering a sparse network such that the fitting error to the data is below a
fixed tolerance.

We formulate the following problem of recovering a sparse network. Given a fixed tolerance, we seek to recover a network such that its $\rms $ is below this tolerance and none of the networks with the same set of vertices as our network and a subset of its edges has a $\rms $ below the tolerance.

\begin{problem}\label{pr2}
Let $V$ be a finite set of vertices, let $E$ be a set of edges on $V$, let $0<|u|_{\rm min}\le |u|_{\rm max}$ be positive values, let $m\in \NN^*$, let $\bm u,\bm s\in \Cc(V,\CC^m)$, respectively $\bm u,\bm s\in  \Cc(V,\RR^m)$, such that $|u|_{\rm min}\leq |u_j| \leq |u|_{\rm max}$ for all $j=1,\ldots,m$ and let $\tol >0$ be a tolerance. Determine an AC network $\Gamma=(V,a)$, respectively a DC network $\Gamma=(V,c)$,  with set of edges $E(\Gamma)\subseteq E$ such that $\rms (\Gamma,\bm u,\bm s)\leq \tol $ and $\Gamma$ is ``minimal'' in the following sense:
Given any electrical network $\Gamma' =(V,a')$, respectively $\Gamma'=(V,c')$, with edge set $E(\Gamma')$, we have that
\begin{enumerate}
\item If $E(\Gamma')= E(\Gamma), \text{then} \, \, \, \rms (\Gamma',\bm u,\bm s)\geq \rms (\Gamma,\bm u,\bm s)$.
\item If $E(\Gamma')\subsetneq E(\Gamma),\text{then} \, \, \, \rms (\Gamma',\bm u,\bm s) > \tol $.
\end{enumerate}
\end{problem}

\begin{remark}
{\rm In particular, we have $\rms (\Gamma',\bm u,\bm s) > \tol  \geq \rms (\Gamma,\bm u,\bm s)$ if $E(\Gamma')\subsetneq E(\Gamma)$. However, notice that there could be an electrical network $E(\Gamma') =(V,a')$ such that
$\rms (\Gamma',\bm u,\bm s)< \rms (\Gamma,\bm u,\bm s)$ and $E(\Gamma')\not\subset E(\Gamma)$.}
\end{remark}

A first naive idea to solve Problem \ref{pr2} could be to recover a network $\Gamma$ solving Problem \ref{fijo} with the set of edges $E$, {\it i.e.}, $[\Gamma,\rms ]=\network\_recovery (E,\bm u,\bm s)$, and then to remove from $E(\Gamma)$ the edges whose admittance values are close to zero. However, this approach presents some problematic issues.

First, let $\Gamma^*$ be a solution to Problem \ref{pr2}, with set of edges $E(\Gamma^*)$. As we will see in Section \ref{s5}, if $E(\Gamma^*)$ is much smaller than $E$, then usually the condition number of the operator associated with the network recovery Problem \ref{fijo} using set $E$, $\kappa(\mathcal{M}_{E,\bm u})$, is much bigger than the condition number of the operator associated with the network recovery Problem \ref{fijo} using set $E(\Gamma^*)$, $\kappa(\mathcal{M}_{E(\Gamma^*),\bm u})$. Then, the values of the admittance of $\Gamma$ at the edges of $E\setminus E(\Gamma^*)$ are usually far from zero. For instance, in Example \ref{ejemplo}, solving the Problem \ref{fijo} with the edge set $E(K_V)=\left\{e_{xy},e_{xz},e_{yz}\right\}$, we can get solutions with zero error such that the values of the admittance at all edges of $E(K_V)$ are far from zero, despite the fact that there are solutions with zero error and less than three edges.

Second, in $\Gamma$ could exist edges with admittance close to zero whose removal would lead to a network which would not correctly fit the data (see Example \ref{ej3}). In order to avoid those problems, the algorithm that we propose in Section \ref{section_algor} algorithm uses techniques of spectral sparsification of networks in order to remove edges from networks.

\section{Spectral network sparsification}\label{sec33}

Spielman and Teng introduced the notion of spectral sparsification of a real
weighted graph, (that is, a DC network), in \cite{sparse_original}, (see also the subsequent papers \cite{sparse_resumen,boletin,sparse} on this topic).

\begin{definition}[\cite{boletin}]\label{esp}
{\rm For $\epsilon > 0$, we say that a DC network $\Gamma'=(V,c')$ with energy $\E'$ is an $\epsilon$-approximation
of a DC network $\Gamma=(V,c)$ with energy $\E$ if for all $u
\in \Cc(V)$:
\begin{equation} \label{dispersa}
	\frac{1}{1+\epsilon}\E(u,u) \leq \E'(u,u) \leq(1+\epsilon) \E(u,u).
\end{equation}}
\end{definition}

In \cite{S24} the definition of the energy was extended to AC networks. Given an AC network $\Gamma = (V, a)$, its energy is a complex bilinear form, so we can not apply Definition \ref{esp} to it. Nevertheless, we can apply this definition separately to its conductance and susceptance networks, $\Gamma_c = (V, c)$ and $\Gamma_b = (V, b)$. We introduce the following definition.

\begin{definition}
{\rm For $\epsilon > 0$, we say that an AC network $\Gamma'$ is an $\epsilon$-approximation of an AC network
$\Gamma$ if the conductance and susceptance networks of $\Gamma'$ are
$\epsilon$-approximations of the conductance and susceptance networks of
$\Gamma$, respectively.}
\end{definition}

In \cite{sparse_resumen} there is a procedure (Algorithm \ref{sparsify}) that can
be used to construct a sparse approximation $\Gamma'=(V,c')$ of a DC network $\Gamma=(V,c)$ such that $E(\Gamma')\subseteq E(\Gamma)$.
\begin{algorithm}\label{spar}
\caption{$\Gamma'= \Sparsify (\Gamma,\epsilon)$.}
\label{sparsify}
\begin{algorithmic}		
\State{Set $t=8|V|\cdot log(|V|)/\epsilon^2$, $c'=0$, and $E(\Gamma')=\emptyset$.}		
\For{each edge $e_{xy}\in E(\Gamma)$}		
\State{Assign to edge $e_{xy}$ a probability $p_{xy}$ proportional to $c(x,y)r^{e}(x,y)$.}		
\EndFor	
\State{Take $t$ samples independently with replacement from $E(\Gamma)$, and each time the edge 
	$e_{xy}$ is sampled, increase the value of $c'(x,y)$ and $c'(y,x)$ by $c(x,y)/tp_{xy}$, (and as a consequence, add the edge $e_{xy}$ to $E(\Gamma')$ if it is the first time that $e_{xy}$ is sampled).}
\end{algorithmic}
\end{algorithm}

Algorithm \ref{sparsify} is not guaranteed to produce an $\epsilon$-approximation of $\Gamma$,
but in \cite{sparse_resumen} there is a proof of the following theorem:
\begin{theorem}[Batson, Spielman, Srivastava, Teng]\label{prob}
Let $\Gamma$ be a DC network, let $\epsilon \in \RR$, $0 <
\epsilon \leq 1$ and let $\Gamma'= \Sparsify (\Gamma,\epsilon)$. Then $\Gamma'$ is an $\epsilon$-approximation of $\Gamma$ with
probability at least $1/2$.
\end{theorem}

Note that the edges that are removed in  Algorithm \ref{sparsify} are the ones that are never sampled. The choice of the sampling probabilities in the algorithm has an interesting physical interpretation. By (\ref{eqsumadc}), for each edge $e_{xy}\in E$, the probability of sampling it, $p_{xy}$, is proportional to the dimensionless ratio
$$
c(x,y)r^{e}(x,y) = \frac{c(x,y)}{c^{e}(x,y)} = \frac{c(x,y)}{c(x,y) +
c_{\setminus e_{xy}}^{e}(x,y)}\in (0, 1],
$$
where $c_{\setminus e_{xy}}^{e}(x,y)\geq 0$ is the effective conductance between $x$ and $y$ in the network $\Gamma\setminus e_{xy}$ obtained from $\Gamma$ by removing the edge $e_{xy}$, which is equal to the contribution to the effective conductance $c^e(x,y)$ of the paths between $x$ and $y$ through the rest of vertices of the network $\Gamma$. 

The probability of choosing edge $e_{xy}$ in each sample is high if the contribution of the value $c(x,y)$ to the effective conductance between $x$ and $y$ in $\Gamma$, $c^e(x,y)$, is very important, that is, if $c(x,y)$ is big compared to $c_{\setminus e_{xy}}^{e}(x,y)$. The limit case is that in which the removal of $e_{xy}$ isolates $x$ and $y$, in which $c_{\setminus e_{xy}}^{e}(x,y)= 0$, so $c(x,y)r^{e}(x,y)=1$; and thus the probability of keeping $e_{xy}$ in $\Gamma'$ is the highest possible.

The probability $p_{xy}$ is low in the opposite case, that is, when $c(x,y)\ll c_{\setminus e_{xy}}^{e}(x,y)$. In that case removal of edge $e_{xy}$ does not significantly affect the effective conductance between $x$ and $y$, because $c_{\setminus e_{xy}}^{e}(x,y) \approx c^{e}(x,y)$. In that case, the probability of keeping $e_{xy}$ in $\Gamma'$ is low.

\begin{example}\label{ejemplo}
{\rm	We consider the DC network in Figure \ref{fisico_grafica}, with conductance $c$ whose values at each edge are given
in Table \ref{tabla_c_ej}. In the same table we indicate the sampling probabilities
calculated for each edge in Algorithm \ref{sparsify}.

\begin{figure}[h!]
	\centering
%\sidecaption
	\includegraphics[width=92mm]{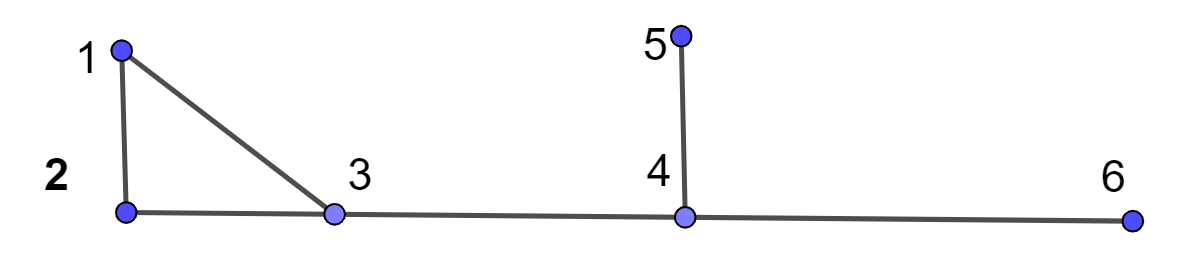}
	\caption{Topology of the network in Example \ref{ejemplo}.}
	\label{fisico_grafica}
\end{figure}

\begin{table}[h!]
%	\centering
	\caption{Conductances and sampling probabilities of the edges in the network.}
	\label{tabla_c_ej}
	\begin{tabular}{p{4.5cm} p{1.5cm} p{1.5cm} p{1.5cm} p{1.5cm} p{1.5cm} p{1.5cm}}
		\hline
		& \multicolumn{6}{c}{Edges} \\
		\cline{2-7}
		& $\{x_1,x_2\}$ & $\{x_1,x_3\}$  & $\{x_2,x_3\}$ & $\{x_3,x_4\}$  & $\{x_4,x_5\}$ & $\{x_4,x_6\}$ \\
		\hline
		Conductance      & $0.5797$    & $75.980$   & $75.980$    & $0.4698$  & $94.599$    & $79.909$\\
		Conductance times effective resistance         & $0.0150$    &  $0.9925$   & $0.9925$   & $1$  & $1$    &  $1$\\
		Sampling probability         & $0.0030$    & $0.1985$   & $0.1985$    & $0.2$  & $0.2$    &  $0.2$\\
		\hline
	\end{tabular}
\end{table}

The values of the conductance at edges $\{x_1,x_2\}$ and $\{x_3,x_4\}$ are two orders of magnitude smaller
than at the rest of edges, nevertheless, edge $\{x_3,x_4\}$ is
crucial because it is the only connection between vertices $x_3$ and $x_4$, so it has a
great probability of being sampled to form any sparse approximation of the network,
while edge $\{x_1,x_2\}$ is expendable, because its conductance $c(x_1,x_2)$ is equal to only $1.5\%$ of the
effective conductance between nodes $x_1$ and $x_2$. The current can flow from node $x_1$ to node
$x_2$ with much greater ease by edges $\{x_1,x_3\}$ and $\{x_2,x_3\}$, so the sampling
probability of $\{x_1,x_2\}$ is close to zero.
}	
\end{example}

In the case of an AC network $\Gamma=(V,a)$, we will denote by $\Gamma'= \Sparsify (\Gamma,\epsilon)$ an AC network constructed following the next procedure. First, we apply Algorithm \ref{sparsify} separately to the conductance and susceptance networks, $\Gamma_c$ and $\Gamma_b$, of $\Gamma$; getting $\Gamma_c'=(V,c')$ and $\Gamma_b'=(V,b')$, respectively. Then, $\Gamma'=(V,a')=(V,c'-ib')$. Note that $E(\Gamma')=E(\Gamma_c') \cup E(\Gamma_b')$, so in the sparsification process of an AC network we remove the edges erased in the sparsification of both $\Gamma_c$ and $\Gamma_b$.

By definition, $\Gamma'$ will be an $\epsilon$-approximation of $\Gamma$ iff $\Gamma_c'$ is an $\epsilon$-approximation of $\Gamma_c$ and $\Gamma_b'$ is an $\epsilon$-approximation of $\Gamma_b$. By Theorem \ref{prob}, if $0 <
\epsilon \leq 1$, then $\Gamma'$ is an $\epsilon$-approximation of $\Gamma$ with probability at least $\frac{1}{4}$.

\section{Sparsification of recovered electrical networks}\label{sec34}

Once we have a procedure to sparsify any network, the key observation to develop our
algorithm of sparse network recovery is that spectral sparsification is
guaranteed to preserve the fitting properties of a network to some extent.

We introduce some notation for the main result of this work in the case of DC networks. Let $\Gamma=(V,c)$ be a DC network with Laplacian $\L$ and let $\bm u,\bm s\in  \Cc(V,\RR^m)$, such that $0<|u|_{\rm min}\leq |u_j| \leq |u|_{\rm max}$ for all $j=1,\ldots,m$. For each $j=1,...,m$, we denote  as $\Qq_{\Gamma, u_j}$ the linear operator on $\Cc(V)$ defined for each $v\in \Cc(V)$ as $\Qq_{\Gamma, u_j}(v)=u_j\L(u_jv)$. Clearly, $\Qq_{\Gamma, u_j}$ is a self-adjoint and positive semidefinite operator. Now, we denote  as $\Qq_{\Gamma, \bm u}$ the linear operator on $\Cc(V,\RR^m)$ defined for each $(v_1,...,v_m)\in \Cc(V,\RR^m)$ as $\Qq_{\Gamma, \bm u}(v_1,...,v_m)=\left(\Qq_{\Gamma, u_1}(v_1),...,\Qq_{\Gamma, u_m}(v_m)\right)$. For any ${\bm v}=(v_1,...,v_m), {\bm w}=(w_1,...,w_m) \in \Cc(V,\RR^m)$, we have that
$$\begin{array}{rl}\displaystyle
\langle \Qq_{\Gamma, \bm u}({\bm v}), {\bm w}\rangle=& \sum_{j=1}^{m}{\int_V \L(u_jv_j) u_jw_jdx}= \sum_{j=1}^{m}{\int_V u_jv_j \L(u_jw_j) dx}\\[3ex]
=& \langle {\bm v},\Qq_{\Gamma, \bm u}({\bm w})\rangle,
\end{array}$$
so the operator $\Qq_{\Gamma, \bm u}$ is self-adjoint. From the same expression, when  ${\bm v}= {\bm w}$, we can see that $\Qq_{\Gamma, \bm u}$ is also positive semidefinite. We denote by $\textbf{1}_m=(\chi_{_V},...,\chi_{_V})\in \Cc(V,\RR^m)$ the vector function whose components are equal to $1$ at all vertices. Note that $\Qq_{\Gamma, \bm u}(\textbf{1}_m)=\mathcal{M}_{E(\Gamma),\bm u}(c)$.

For each $j=1,...,m$, we denote  as $u_j^{-1}\in \Cc(V)$ the function defined for each $x\in V$ as $u_j^{-1}(x)=1/u_j(x)$. We also define the vector function $$
\phi_{\bm u} = \left(\frac{\int_V{u_1^{-1}  dx}}{||u_1^{-1}||^2} u_1^{-1},...,\frac{\int_V{u_m^{-1} dx}}{||u_m^{-1}||^2} u_m^{-1}\right)\in \Cc(V,\RR^m).
$$

Our main original result consists in the following upper bound for the $\rms $ of any
sparse approximation in a fixed data set.

\begin{theorem}[Main theorem] \label{cota}
Given a DC network $\Gamma=(V,c)$, positive values $0<|u|_{\rm min}\le |u|_{\rm max}$, $m\in \NN^*$ and $\bm u,\bm s\in  \Cc(V,\RR^m)$ such that $|u|_{\rm min}\leq |u_j| \leq |u|_{\rm max}$ for all $j=1,\ldots,m$; if $\Gamma'$ is an $\epsilon$-approximation of $\Gamma$, then:
\[
\rms (\Gamma',\bm u,\bm s) \leq \rms (\Gamma,\bm u,\bm s) + \epsilon
\frac{\|\Qq_{\Gamma, \bm u}\|_2\cdot\| \textbf{1}_m-\phi_{\bm u}\|}{\sqrt{m|V|}}.
\]
\end{theorem}

\begin{proof}
Let $\Gamma'=(V,c')$ be an $\epsilon$-approximation of $\Gamma$, with Laplacian $\L'$. As $\Qq_{\Gamma', \bm u}(\textbf{1}_m)=\mathcal{M}_{E(\Gamma'),\bm u}(c')$, by (\ref{eqrms}) we get that
\begin{equation}\notag
\begin{split}
	\rms (\Gamma',\bm u,\bm s)  & = \frac{1}{\sqrt{m|V|}}\|\Qq_{\Gamma', \bm u}(\textbf{1}_m) - {\bm s}\|  \\
	& \leq \frac{1}{\sqrt{m|V|}}\|\Qq_{\Gamma, \bm u}(\textbf{1}_m) - {\bm s}\| +
	\frac{1}{\sqrt{m|V|}}\|\Qq_{\Gamma', \bm u}(\textbf{1}_m)-\Qq_{\Gamma, \bm u}(\textbf{1}_m)\| \\
	& = \rms (\Gamma,\bm u,\bm s)  +
	\frac{1}{\sqrt{m|V|}}\|\left(\Qq_{\Gamma', \bm u}-\Qq_{\Gamma, \bm u}\right)(\textbf{1}_m)\|.
\end{split}
\end{equation}

The constant functions belong to the
null space of any Laplacian, therefore $\phi_{\bm u}$ belongs to the null space of
$\Qq_{\Gamma', \bm u}-\Qq_{\Gamma, \bm u}$, and thus
\begin{equation} \label{truco}
\begin{split}
	\|\left(\Qq_{\Gamma', \bm u}-\Qq_{\Gamma, \bm u}\right)(\textbf{1}_m)\|  &=
	\|\left(\Qq_{\Gamma', \bm u}-\Qq_{\Gamma, \bm u}\right)
	\left(\textbf{1}_m-\phi_{\bm u}\right)\|  \\
	&\leq \|\Qq_{\Gamma', \bm u}-\Qq_{\Gamma, \bm u}\|_2 \cdot
	\left\|\textbf{1}_m-\phi_{\bm u}\right\|.
\end{split}
\end{equation}

In order to complete the proof, it is enough to show that
$\|\Qq_{\Gamma', \bm u}-\Qq_{\Gamma, \bm u}\|_2 \leq \epsilon
\|\Qq_{\Gamma, \bm u}\|_2$. Now, $\Qq_{\Gamma', \bm u}-\Qq_{\Gamma, \bm u}$ is a self-adjoint operator, so its spectral norm is equal to the maximum of the absolute value of its eigenvalues. It is straightforward to prove that its eigenvalues are those of the operators $\Qq_{\Gamma', u_j}-\Qq_{\Gamma, u_j}$ for all $1 \leq j \leq m$. In particular, there is at least one index $k$ such that
\[
\|\Qq_{\Gamma', \bm u}-\Qq_{\Gamma, \bm u}\|_2=\|\Qq_{\Gamma', u_k}-\Qq_{\Gamma, u_k}\|_2.
\]

Let $u_0\in \Cc(V)$ be a eigenvector of $\Qq_{\Gamma', u_k}-\Qq_{\Gamma, u_k}$
corresponding with its eigenvalue that has the largest absolute value such that $\|u_0\|=1$. Then
\[
\|\Qq_{\Gamma', \bm u}-\Qq_{\Gamma, \bm u}\|_2= |\langle u_0,
\left(\Qq_{\Gamma', u_k}-\Qq_{\Gamma, u_k}\right) (u_0)\rangle|.
\]
Now, denoting as $\E$ and $\E'$ the energy of $\Gamma$ and $\Gamma'$, respectively, we have that
\begin{equation}\notag
\begin{split}
	\langle u_0,\left(\Qq_{\Gamma', u_k}-\Qq_{\Gamma, u_k}\right) (u_0)\rangle  & = \langle u_0,u_k(\L'-\L)(u_ku_0)\rangle \\
	& = \int_V u_0 u_k(\L'-\L)(u_0 u_k) dx\\
	& = \E'(u_0 u_k,u_0 u_k)-\E(u_0 u_k,u_0 u_k).
\end{split}
\end{equation}

Evaluating (\ref{dispersa}) at $u=u_0u_k\in \Cc(V)$ we get
\begin{equation}\label{dem}
\frac{1}{1+\epsilon}\E(u_0u_k,u_0u_k) \leq \E'(u_0u_k,u_0u_k) \leq(1+\epsilon) \E(u_0u_k,u_0u_k).
\end{equation}
From the right inequality of (\ref{dem}), we have:
\[
\E'(u_0 u_k,u_0 u_k)-\E(u_0 u_k,u_0 u_k) \leq \epsilon \E(u_0 u_k,u_0 u_k),
\]
and from the left inequality in (\ref{dem}):
\[
-\left(\E'(u_0 u_k,u_0 u_k)-\E(u_0 u_k,u_0 u_k)\right) \leq \frac{\epsilon}{1+\epsilon} \E(u_0 u_k,u_0 u_k) < \epsilon \E(u_0 u_k,u_0 u_k).
\]

Joining the last two inequalities, we get
\[
\|\Qq_{\Gamma', \bm u}-\Qq_{\Gamma, \bm u}\|_2 \leq \epsilon \E(u_0 u_k,u_0 u_k).
\]
Moreover,
\begin{equation}\notag
\begin{split}
	\E(u_0 u_k,u_0 u_k) & = \int_V u_0 u_k \L(u_0 u_k) dx \\
	& = \langle u_0,u_k \L(u_ku_0)\rangle\\
	& = \langle u_0,\Qq_{\Gamma, u_k}(u_0)\rangle.
\end{split}
\end{equation}

By the Courant-Fisher theorem \cite{sparse}, the evaluation of the quadratic form
associated with the operator $\Qq_{\Gamma, u_k}$ at any unit vector is less or equal than its
maximum eigenvalue, which is equal to the norm of this operator because it is positive
semidefinite, so:
\begin{equation}\label{eqpruebamain}
\|\Qq_{\Gamma', \bm u}-\Qq_{\Gamma, \bm u}\|_2 \leq \epsilon \|
\Qq_{\Gamma, u_k} \|_2 \leq  \epsilon \|\Qq_{\Gamma, \bm u}\|_2,
\end{equation}
where the last inequality holds because the eigenvalues of $\Qq_{\Gamma, \bm u}$ are those of the operators $\Qq_{\Gamma, u_j}$ for all $1 \leq j \leq m$.
\end{proof}

\begin{remark}
{\rm In (\ref{truco}), it is possible to introduce any vector function equal to 
$\left(\lambda_1 u_1^{-1},...,\lambda_m u_m^{-1}\right)\in \Cc(V,\RR^m)$, for any choice of $\lambda_1,...,\lambda_m\in \RR$, because all of them belong to the null space of $\Qq_{\Gamma', \bm u}-\Qq_{\Gamma, \bm u}$. Among them, we choose to
introduce the term $\phi_{\bm u}$ because it is the choice that minimizes
$\rho(\lambda_1,..,\lambda_m) \equiv \\\left\|\textbf{1}_m-\left(\lambda_1 u_1^{-1},...,\lambda_m u_m^{-1}\right)\right\|^2=\sum_{j=1}^{m}{\int_V \left(\chi_{_V}- \lambda_j u_j^{-1}\right)^2 dx}$.
Effectively, if we look for a minimum of $\rho$, its partial derivatives must be
zero:
\[
\frac{\partial \rho(\lambda_1,..,\lambda_m)}{\partial \lambda_k} = 2
\int_V{\left(-u_k^{-1} + \lambda_k\left(u_k^{-1}\right)^2\right) dx} = 0,
\]
so the only possible choice for the $\lambda_k$ is:
\[
\lambda_k = \frac{\int_V{u_k^{-1} dx}}{||u_k^{-1}||^2},
\]
which is indeed the global minimum of $\rho$ because its Hessian is positive definite: if $j\neq k$, then $\frac{\partial^2 \rho(\lambda_1,..,\lambda_m)}{\partial \lambda_k \partial \lambda_j}=0$, and for each $k=1,...,m$,
\[
\frac{\partial^2 \rho(\lambda_1,..,\lambda_m)}{\partial \lambda_k^2} = 2\int_V{\left(u_k^{-1}\right)^2 dx}> 0.
\]}
\end{remark}

\begin{remark}
{\rm The term $\frac{\|\Qq_{\Gamma, \bm u}\|_2\cdot\| \textbf{1}_m-\phi_{\bm u}\|}{\sqrt{m|V|}}$ that appears in the upper bound given by the
last theorem depends only on the electrical network $\Gamma$ and the first element of the data pair $(\bm u,\bm s)$, so if we want a sparse approximation of $\Gamma$ such that the $\rms $ of
the approximation does not increase (with respect to $\rms (\Gamma,\bm u,\bm s)$) in the
sparsification procedure above a fixed value, there is an $\epsilon > 0$ such
that any $\epsilon$-approximation of $\Gamma$  is guaranteed to meet this
requirement.}
\end{remark}

In the case of an AC electrical network $\Gamma$, we denote by $\Gamma_c=(V,c)$ and $\Gamma_b=(V,b)$ the conductance and susceptance networks of $\Gamma$, whose respective Laplacians are $\L_c$ and $\L_b$. Additionally, we define the function
\begin{equation}
\resizebox{\hsize}{!}{$\begin{split}
	\Delta (\Gamma,\bm u) & = \left( \left(
	\sqrt{m|V|}(\|\Qq_{\Gamma_c, \Re(\bm u)}\|_2+\|\Qq_{\Gamma_c, \Im(\bm u)}\|_2) +
	\|\L_b\|_2\left(\|\Re(\bm u)\|_{\infty}
	\|\Im(\bm u)\|_2+\|\Im(\bm u)\|_{\infty} \|\Re(\bm u)\|_2\right)\right)^2 +
	\right. \\
	&  \left. + \left( \sqrt{m|V|}(\|
	\Qq_{\Gamma_b, \Re(\bm u)}\|_2
	+\|\Qq_{\Gamma_b, \Im(\bm u)}\|_2) +
	\|\L_c\|_2(\|\Im(\bm u)\|_{\infty}
	\|\Re(\bm u)\|_2+\|\Re(\bm u)\|_{\infty} \|\Im(\bm u)\|_2)\right)^2\right)^{1/2}, \notag
\end{split}$}
\end{equation}
and we have the following result, which is the AC analogous to Theorem \ref{cota}:

\begin{theorem} \label{cota_ac}
Given an AC network $\Gamma=(V,a)$, $m\in \NN^*$ and $\bm u,\bm s\in  \Cc(V,\CC^m)$, if $\Gamma'$ is an $\epsilon$-approximation of $\Gamma$, then:
\[
\rms (\Gamma',\bm u,\bm s) \leq \rms (\Gamma,\bm u,\bm s) + \epsilon \frac{\Delta (\Gamma,\bm u)}{\sqrt{2m|V|}},
\]
where $\Delta (\Gamma,\bm u)$ is a function which depends only on $\Gamma$ and $\bm u$.
\end{theorem}

\begin{proof}
First, for every $u\in \Cc(V,\CC)$, we have that
\begin{equation}\notag
u\overline{\L(u)}=u\L^*(\overline{u})=(\Re(u)+i\Im(u))(\L_c+i\L_b)((\Re(u)-i\Im(u))),
\end{equation}
and therefore
\begin{equation} \label{eqprueba1}
\Re(u\overline{\L(u)})=\Re(u)\L_c(\Re(u))+  \Re(u)\L_b(\Im(u))+  \Im(u)\L_b(\Re(u))+  \Im(u)\L_c(\Im(u))  ,
\end{equation}
and also
\begin{equation} \label{eqprueba2}
\Im(u\overline{\L(u)})=\Im(u)\L_c(\Re(u))+  \Im(u)\L_b(\Im(u))+  \Re(u)\L_b(\Re(u))+  \Re(u)\L_c(\Im(u)).
\end{equation}

Now, let $\Gamma'=(V,a')$ be an $\epsilon$-approximation of $\Gamma$, with Laplacian $\L'$. Let $\Gamma_c'=(V,c')$ and $\Gamma_b'=(V,b')$ be the conductance and susceptance networks of $\Gamma'$, whose respective Laplacians are $\L_c'$ and $\L_b'$. By (\ref{eqrms}) we get that
\begin{equation}\notag
\begin{split}
	\rms (\Gamma',\bm u,\bm s)  & = \frac{1}{\sqrt{m|V|}}\|\mathcal{M}_{E(\Gamma'),\bm u}(c',b') - \mathcal{S}(\bm s)\|  \\
	& \leq \frac{1}{\sqrt{m|V|}}\|\mathcal{M}_{E(\Gamma),\bm u}(c,b) - \mathcal{S}(\bm s)\| \\
	& +
	\frac{1}{\sqrt{m|V|}}\|\mathcal{M}_{E(\Gamma'),\bm u}(c',b')-\mathcal{M}_{E(\Gamma),\bm u}(c,b)\| \\
	& = \rms (\Gamma,\bm u,\bm s)  +
	\frac{1}{\sqrt{m|V|}}\|\mathcal{M}_{E(\Gamma'),\bm u}(c',b')-\mathcal{M}_{E(\Gamma),\bm u}(c,b)\|.
\end{split}
\end{equation}

The square of the norm in the last equation is equal to
\begin{align}\label{eqprueba3}
% & \|\mathcal{M}_{E(\Gamma'),\bm u}(c',b')-\mathcal{M}_{E(\Gamma),\bm u}(c,b)\|^2 \notag \\ 
% & = \|\big(\Re(u_1\overline{\left(\L'-\L\right)(u_1)}), ... ,\Re(u_m\overline{\left(\L'-\L\right)(u_m)})\big)\|^2 \notag \\ 
% &+ \|\big(\Im(u_1\overline{\left(\L'-\L\right)(u_1)}), ... ,\Im(u_m\overline{\left(\L'-\L\right)(u_m)})\big)\|^2.
\|\mathcal{M}_{E(\Gamma'),\bm u}(c',b')-\mathcal{M}_{E(\Gamma),\bm u}(c,b)\|^2 &= \|\big(\Re(u_1\overline{\left(\L'-\L\right)(u_1)}), ... ,\Re(u_m\overline{\left(\L'-\L\right)(u_m)})\big)\|^2 \notag \\ 
 &+ \|\big(\Im(u_1\overline{\left(\L'-\L\right)(u_1)}), ... ,\Im(u_m\overline{\left(\L'-\L\right)(u_m)})\big)\|^2.
\end{align}

By equations (\ref{eqprueba1}), (\ref{eqprueba2}) and (\ref{eqprueba3}), we can write
\begin{equation}  \label{eqprueba4}
\begin{split}
	\|\mathcal{M}_{E(\Gamma'),\bm u}(c',b')-\mathcal{M}_{E(\Gamma),\bm u}(c,b)\|^2 &  =
	\| \beta_1+\beta_2-\beta_3+\beta_4\|^2\\
	& +\| \beta_5+\beta_6+\beta_7-\beta_8\|^2 \\
	& \leq \left( \|\beta_1\|+\|\beta_2\|+\|\beta_3\|+\|\beta_4\|\right)^2\\
	& + \left( \|\beta_5\|+\|\beta_6\|+\|\beta_7\|+\|\beta_8\|\right)^2,
\end{split}
\end{equation}
where:

\vspace{0.3cm}
\begin{tabular}{ll} 
	$\beta_1=\left(\Pp_{c'}-\Pp_{c}\right)\left(\Re(\bm u),\Re(\bm u)\right),$
	&\hspace{.8cm}
	$\beta_2=\left(\Pp_{b'}-\Pp_{b}\right)\left(\Re(\bm u),\Im(\bm u)\right),$\\
	$\beta_3=\left(\Pp_{b'}-\Pp_{b}\right)\left(\Im(\bm u),\Re(\bm u)\right),
	$&\hspace{.8cm}
	$\beta_4=\left(\Pp_{c'}-\Pp_{c}\right)\left(\Im(\bm u),\Im(\bm u)\right),$\\
	$\beta_5=\left(\Pp_{c'}-\Pp_{c}\right)\left(\Im(\bm u),\Re(\bm u)\right),$
	&\hspace{.8cm}
	$\beta_6=\left(\Pp_{b'}-\Pp_{b}\right)\left(\Im(\bm u),\Im(\bm u)\right),$\\
	$\beta_7=\left(\Pp_{b'}-\Pp_{b}\right)\left(\Re(\bm u),\Re(\bm u)\right),$
	&\hspace{.8cm}
	$\beta_8=\left(\Pp_{c'}-\Pp_{c}\right)\left(\Re(\bm u),\Im(\bm u)\right).$
\end{tabular}
\vspace{0.3cm}

Then, on one hand, $\Pp_{c}\left(\Re(\bm u),\Re(\bm u)\right)=\Qq_{\Gamma_c, \Re(\bm u)}(\textbf{1}_m)$ and $\Pp_{c'}\left(\Re(\bm u),\Re(\bm u)\right)=\Qq_{\Gamma_c', \Re(\bm u)}(\textbf{1}_m)$, so by (\ref{eqpruebamain}), we get 
\begin{equation} \notag
\begin{split}
	\|\beta_1\|  &\leq \|\Qq_{\Gamma_c', \Re(\bm u)}-\Qq_{\Gamma_c, \Re(\bm u)}\|_2 \cdot \left\|\textbf{1}_m\right\| =\sqrt{m|V|}\|\Qq_{\Gamma_c', \Re(\bm u)}-\Qq_{\Gamma_c, \Re(\bm u)}\|_2\\
	&\leq \epsilon \sqrt{m|V|}\|\Qq_{\Gamma_c, \Re(\bm u)}\|_2.
\end{split}
\end{equation}

Reasoning analogously, we also obtain that $\|\beta_4\| \leq \epsilon \sqrt{m|V|}\|\Qq_{\Gamma_c, \Im(\bm u)}\|_2$, $\|\beta_6\| \leq \epsilon \sqrt{m|V|}\|\Qq_{\Gamma_b, \Im(\bm u)}\|_2$ and $\|\beta_7\| \leq \epsilon \sqrt{m|V|}\|\Qq_{\Gamma_b, \Re(\bm u)}\|_2$.

On the other hand, $\L_c=\Qq_{\Gamma_c, \chi_{_V}}$, $\L_c'=\Qq_{\Gamma_c', \chi_{_V}}$, $\L_b=\Qq_{\Gamma_b, \chi_{_V}}$ and $\L_b'=\Qq_{\Gamma_b', \chi_{_V}}$, so, by the proof of Theorem \ref{cota}, we have that $\|\L_c'-\L_c\|_2\leq \epsilon \|\L_c\|_2$ and $\|\L_b'-\L_b\|_2\leq \epsilon \|\L_b\|_2$. As a consequence, for any $({\bm v},{\bm w})=\left((v_1,...,v_m),(w_1,...,w_m)\right)\in \Cc(V,\RR^m)\times\Cc(V,\RR^m)$ it is satisfied that 
\begin{equation} \notag
\begin{split}
	\|(\Pp_{c'}-\Pp_{c})\left({\bm v},{\bm w}\right)\|^2  &= \sum_{j=1}^{m}{\int_V v_j^2(\L_c'-\L_c)(w_j)^2 dx}\\
	&\leq \|{\bm v}\|_{\infty}^2 \sum_{j=1}^{m}{\|(\L_c'-\L_c)(w_j)\|_2^2}\\
	&\leq \|{\bm v}\|_{\infty}^2 \cdot \|(\L_c'-\L_c)\|_2^2 \cdot \|{\bm w}\|_2^2\\
	&\leq \epsilon^2 \|\L_c\|_2^2 \cdot \|{\bm v}\|_{\infty}^2 \cdot  \|{\bm w}\|_2^2.
\end{split}
\end{equation}

And, analogously, it is satisfied that $\|(\Pp_{b'}-\Pp_{b})\left({\bm v},{\bm w}\right)\|\leq \epsilon \|\L_b\|_2 \cdot \|{\bm v}\|_{\infty} \cdot  \|{\bm w}\|_2$. Applying this result, we get the following bounds: %\\$\|\beta_2\| \leq \|\L_b\|_2 \cdot \|\Re(\bm u)\|_{\infty} \cdot  \|\Im(\bm u)\|_2$, \, \, \, \, \, \, \, \hspace{0.01cm}  $\|\beta_3\| \leq \|\L_b\|_2 \cdot \|\Im(\bm u)\|_{\infty} \cdot  \|\Re(\bm u)\|_2$, \, \, \, \\$\|\beta_5\| \leq \|\L_c\|_2 \cdot \|\Im(\bm u)\|_{\infty} \cdot  \|\Re(\bm u)\|_2$, \, \, and  \, \, $\|\beta_8\| \leq \|\L_c\|_2 \cdot \|\Re(\bm u)\|_{\infty} \cdot  \|\Im(\bm u)\|_2$.

\begin{tabular}{ll} 
	$\|\beta_2\| \leq \|\L_b\|_2 \cdot \|\Re(\bm u)\|_{\infty} \cdot  \|\Im(\bm u)\|_2$,
	&\hspace{.8cm}
	$\|\beta_3\| \leq \|\L_b\|_2 \cdot \|\Im(\bm u)\|_{\infty} \cdot  \|\Re(\bm u)\|_2$,\\
	$\|\beta_5\| \leq \|\L_c\|_2 \cdot \|\Im(\bm u)\|_{\infty} \cdot  \|\Re(\bm u)\|_2$,&\hspace{.8cm}
	$\|\beta_8\| \leq \|\L_c\|_2 \cdot \|\Re(\bm u)\|_{\infty} \cdot  \|\Im(\bm u)\|_2$.
\end{tabular}

Considering in (\ref{eqprueba4}) all the upper bounds obtained for $\|\beta_1\|,...,\|\beta_8\|$, we conclude the proof.
\end{proof}

\section{Algorithm for sparse network recovery}\label{section_algor}

In this section we propose an algorithm to solve the Problem \ref{pr2} of recovering simultaneously the topology and admittance of a sparse electrical network.

If we have an
electrical network $\Gamma =(V,a)$, a data pair  $(\bm u,\bm s)$ and we choose an
$\epsilon$ such that the $\rms $ of any $\epsilon$-approximation of $\Gamma$ in the data pair $(\bm u,\bm s)$ does not
surpass a fixed tolerance $\tol>0$ (by Theorems \ref{cota} and \ref{cota_ac} such a
$\epsilon$ exists if $\rms (\Gamma,\bm u,\bm s)< \tol $), the only guarantee about the number of edges that any $\epsilon$-approximation
obtained by executing $\Sparsify (\Gamma,\epsilon)$ will have is that this number of edges will be less or
equal than the number of edges sampled in Algorithm \ref{sparsify}, $t=8|V|\cdot \log(|V|)/\epsilon^2$. If $\tol$ is small, the number of edges sampled will be large, so most
$\epsilon$-approximations of $\Gamma$ will have the same number of edges as the
original network, thus they will not be useful for our purposes.

In the experimentation we have found that, if we choose an $\epsilon$ such that $\Gamma'=\Sparsify (\Gamma,\epsilon)$ satisfies that $|E(\Gamma')|<|E(\Gamma)|$, and then we solve Problem \ref{fijo} with the set $E(\Gamma')$; (that is, we execute $[\Gamma'', \rms '']=\network\_recovery (E(\Gamma'),\bm u,\bm s)$, determining a network $\Gamma''=(V,a'')$ with set of edges $E(\Gamma'')\subseteq E(\Gamma')$ such that the error $\rms ''=\rms (\Gamma'', \bm u,\bm s)$ is minimum), then we have $\rms''\leq \tol $ in many cases, even in cases in
which $\rms (\Gamma',\bm u,\bm s)> \tol $.

Our approach to solve Problem \ref{pr2} consists in the application of
Algorithm \ref{selfhealing}. The algorithm starts by fitting a network $\Gamma =\network\_recovery (E,\bm u,\bm s)$ solving Problem \ref{fijo} with the data pair  $(\bm u,\bm s)$ and the set of edges $E$. Then, our
goal is to perform a sparsification of $\Gamma$, $\Gamma'= \Sparsify (\Gamma,\epsilon)$, followed by a recovery of another network $\Gamma''$ solving Problem \ref{fijo} with the set of edges of the sparse approximation, $E(\Gamma')$, as in the paragraph above, looking for a
network with less edges than the current one, and with a $\rms $ below the input
tolerance $\tol $.

We do not know which choices of $\epsilon$ can lead us to a network with these
characteristics, so we use a procedure of exploration to let Algorithm \ref{selfhealing} find
suitable values of $\epsilon$. We start the algorithm with an initial input value
of $\epsilon$, and we do an iterative procedure. Each iteration begins by checking
if $\Gamma'= \Sparsify (\Gamma,\epsilon)$ has less edges than the current network
$\Gamma$. If this is not the case, there is a high probability, from Theorem \ref{prob},
of  $\Gamma'$ being a sparse approximation so close to $\Gamma$ that it has the
same set of edges, so we increase the value of $\epsilon$ multiplying it by the
input parameter  $\psi>1$, and we finish the iteration. In this way, in the next
iteration the number of edges sampled in $\Sparsify$ will be lower than in the
previous one, increasing the probability of getting sparse approximations with less
edges than $\Gamma$.

If, on the contrary, $|E(\Gamma')|<|E(\Gamma)|$, we recover the network $\Gamma''$ solving Problem \ref{fijo} with set of edges $E(\Gamma')$, that is, we compute $[\Gamma'', \rms '']=\network\_recovery (E(\Gamma'),\bm u,\bm s)$.

Then, we check if $\rms''\leq\tol$.
If this is the case, we replace the current network $\Gamma$ by this new
$\Gamma''$ and we finish the
iteration, using the new network as input for $\Sparsify$ in the next iteration, repeating the process, in
order to look for networks with less edges than it. If, on the opposite case,
$\rms ''>\tol $, this means we have removed edges that are necessary for the network to correctly fit the data in the sparsification process, because there are no networks with set of edges contained in $E(\Gamma')$
whose $\rms$ is lower or equal than $\tol$. We reject the network $\Gamma''$, we decrease the value of $\epsilon$
dividing it by $\psi$ and we finish the iteration because, from theorems \ref{cota} and \ref{cota_ac},
we know that decreasing the value of $\epsilon$ assures that the $\rms $ on any
subsequent $\epsilon$-approximation of $\Gamma$ will have a smaller upper bound.

It is possible to define different stopping criteria for Algorithm \ref{selfhealing}, such as the
total running time if we are interested on the best network that the algorithm can
find in that period, or a maximum number of consecutive iterations in which the
network has not changed. In the case of distribution networks, it is usual that the
position of the switches is set so that the network topology is a tree, so in this
case, reaching a tree topology can be another stopping criterion. A
detailed discussion about the stopping criteria is left for future work.

\begin{algorithm}
\caption{Sparse network recovery.}
\label{selfhealing}
\begin{algorithmic}		
\State{Inputs: $(\bm u,\bm s)$, $E$, $\epsilon$, $\psi>1$,
	$\tol $, $total\_time$.}		
\State{$[\Gamma, \rms ]:= \network\_recovery (E,\bm u,\bm s)$.}	
\While{$running\_time < total\_time$  }	// (Stopping criteria)	
\State{Set $\Gamma':=\Sparsify (\Gamma,\epsilon)$.}		
\If{$|E(\Gamma')| < |E(\Gamma)|$}		
\State{$[\Gamma'', \rms '']=\network\_recovery (E(\Gamma'),\bm u,\bm s)$.}		
\If{$\rms '' \leq \tol$}	
\State{$\Gamma:=\Gamma''$.}		
\Else	
\State{$\epsilon:=\epsilon/\psi$.}		
\EndIf		
\Else	
\State{$\epsilon=\epsilon\cdot \psi$.}		
\EndIf	
\EndWhile		
\State Return $\Gamma$.		
\end{algorithmic}
\end{algorithm}

\section{Experimental results and discussion}\label{s5}
Lastly, we present some results of the application of Algorithm \ref{selfhealing} to solve examples of Problem \ref{pr2}. The algorithm has been written in MATLAB using Casadi \cite{Ca23}, an open-source software tool that provides a symbolic framework suited for numerical optimization, and the interior-point solver IPOPT \cite{Ip23}, an open-source software package for large-scale nonlinear optimization, for the $\network\_recovery$ process. The tolerance used in IPOPT is equal to $10^{-8}$. In all the examples, we use $\psi=1.5$ as input of Algorithm \ref{selfhealing}.

In each example, we have chosen a network $\Gamma$ and we have sampled a data pair  $(\bm u,\bm s)$ consisting in $m=1000$ pairs of voltage and power injected at $\Gamma$. Each pair $\left(u_j,s_j\right)$ has been computed numerically solving the power flow equations of $\Gamma$, $s_j=u_j\overline{\L(u_j)}$, along with different additional restrictions. Restrictions that are common for all the pairs $\left(u_j,s_j\right)$, such as for instance having $s_j(x)=0$ for a vertex $x\in V$ for all $i=1,...m$, will be highlighted only if they are relevant for the analysis of the results.

In all the examples, we solve Problem \ref{pr2} using $E=E(K_{V})$ as set of edges, that is, we solve the problem without any {\it a priori} information about the topology of the network. We fix a labeling $\left\{x_1,...,x_n\right\}$ on the set of vertices of each network and in the figures we denote each vertex $x_j$ as $j$ for the sake of clarity.

\begin{example}\label{ej3}
{\rm	
We start with the DC network in Example \ref{ejemplo}, (whose conductance is shown in Table \ref{tabla_c_ej}). In the data pair  $(\bm u,\bm s)$, for each $j=1,...,m$ we have different loads, (that is, negative power injected $s_j$) at all nodes, except at node $x_1$, where the power feeding the entire network is generated ($s_j(x_1)\geq0$). We apply Algorithm \ref{selfhealing} with $\tol=10^{-5}$ and an initial value for $\epsilon$ equal to $0.1$ and we obtain the results of Table \ref{fisico_tabla} (we only write the results of iterations in which the number of edges has decreased, and we consider the first network recovery with the set of edges $E(K_{V})$ of the complete graph previous to the iterative process as the iteration number $1$).

\begin{table}[h!]
	\centering{
		\caption{Example \ref{ej3} network fitting results with tolerance $10^{-5}$.\textcolor{white}{$\frac{}{}$}}
		\label{fisico_tabla}
\begin{tabular}{p{3cm} p{1.5cm} p{2.5cm} p{2.5cm} p{1.5cm}}			\hline
			Iteration number &   $|E(\Gamma)|$  &  $\rms(\Gamma,\bm u,\bm s)$  &  $\kappa(\mathcal{M}_{E(\Gamma),\bm u})$  &  $\epsilon$ \\
			\hline
			1 &   15  &   $1.011\cdot 10^{-6}$   &   $1.040\cdot 10^{4}$   &    0.1000\\
			2  &  9 &   $9.415\cdot 10^{-7}$  &   $4.316\cdot 10^{3}$    &   0.1000\\
			3  &   7  &   $4.632\cdot 10^{-7}$   &   $3.021\cdot 10^{3}$   &   0.1000\\
			4  &   6  &   $1.491\cdot 10^{-8}$  &   $2.100\cdot 10^{3}$    &   0.1000\\
			\hline
	\end{tabular}}
\end{table}

The last row in Table \ref{fisico_tabla} corresponds to the final network found by the algorithm, that has the real topology of the network (Figure \ref{fisico_grafica}). The difference between the values at each edge of the conductance obtained and of the real conductance is of the order of magnitude of $10^{-4}$. We obtain almost the real network in just $4$ iterations.

\begin{table}[h!]
	\centering
	\caption{Example \ref{ej3} network fitting results with tolerance $10^{-3}$.\textcolor{white}{$\frac{}{}$}}
	\label{fisico_tabla_sparse}
\begin{tabular}{p{3cm} p{1.5cm} p{2.5cm} p{2.5cm} p{1.5cm}}		\hline
		Iteration number &   $|E(\Gamma)|$  &  $\rms(\Gamma,\bm u,\bm s)$  &  $\kappa(\mathcal{M}_{E(\Gamma),\bm u})$  &  $\epsilon$ \\
		\hline
		1 &   15  &   $1.011\cdot 10^{-6}$   &   $1.040\cdot 10^{4}$   &    0.1000\\
		2  &  11 &   $8.418\cdot 10^{-7}$  &   $6.354\cdot 10^{3}$    &   0.1000\\
		3  &   6  &   $3.580\cdot 10^{-5}$   &   $1.320\cdot 10^{3}$   &   0.1000\\
		5  &   5 &   $3.579\cdot 10^{-5}$  &   $7.028\cdot 10^{2}$    &   0.1500\\
		\hline
	\end{tabular}
\end{table}

If we repeat the experiment with a looser tolerance $\tol=10^{-3}$, we obtain the results in Table \ref{fisico_tabla_sparse}. The final network is a sparse approximation of the original, with the same edges than it except for edge $\{x_1,x_2\}$, that is removed. This removal barely distorts the structure of the network, as we discussed in Example \ref{ejemplo}. It is important to remark that the algorithm never removes the low conductance edge $\{x_3,x_4\}$. That removal would make impossible satisfying the demand of power at nodes $x_4$, $x_5$ and $x_6$, because node $x_1$ is the only node with generated power, so the $\rms$ of the resulting network would be higher than $\tol$.

}
\end{example}

From now on, we will use $\tol=10^{-5}$ for all experiments.

\begin{example} \label{ej}
{\rm	The next experiment uses the CIGRE test AC network named ``Medium voltage distribution network with PV and Wind DER" \cite{ejemplo_cigre}.

\begin{figure}[h!]
	\centering
	\includegraphics[width=160mm]{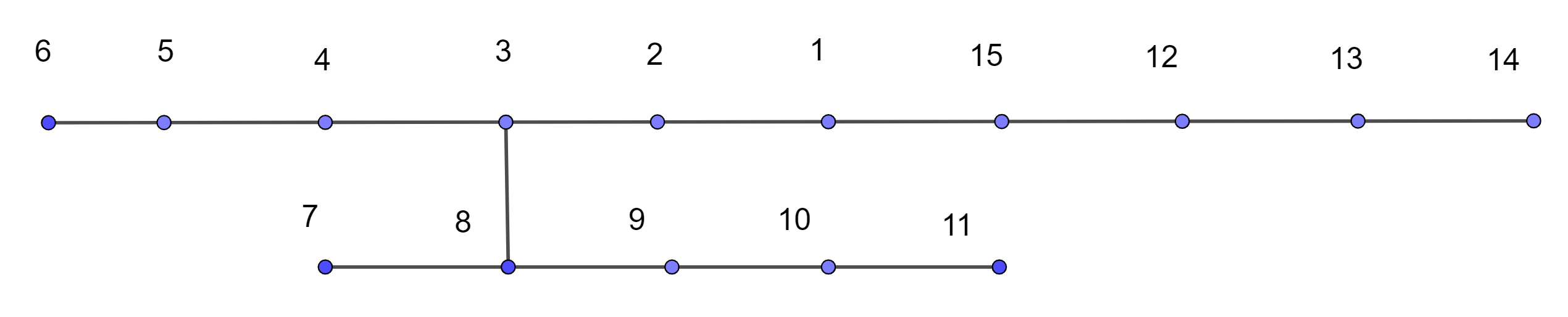}
	\caption{Topology of the CIGRE Network. (Example \ref{ej}).}
	\label{cigre}
\end{figure}

An application of Algorithm \ref{selfhealing} gives the results in Table \ref{cigre_tabla1}.

\begin{table}[h!]
	\centering
	\caption{Example \ref{ej} network fitting results.}
	\label{cigre_tabla1}
\begin{tabular}{p{3cm} p{1.5cm} p{2.5cm} p{2.5cm} p{1.5cm}}		\hline
		Iteration number &   $|E(\Gamma)|$  &  $\rms(\Gamma,\bm u,\bm s)$  &  $\kappa(\mathcal{M}_{E(\Gamma),\bm u})$  &  $\epsilon$ \\
		\hline
		1 &   105  &   $2.905\cdot 10^{-5}$   &   $4.200\cdot 10^{15}$   &    0.3000\\
		2  &  24 &   $8.772\cdot 10^{-6}$  &   $1.492\cdot 10^{15}$    &   0.3000\\
		3  &   18  &   $6.405\cdot 10^{-6}$   &   $1.526\cdot 10^{15}$   &   0.3000\\
		4  &   17  &   $6.397\cdot 10^{-6}$  &   $1.497\cdot 10^{15}$    &   0.3000\\
		6  &  16 &   $4.343\cdot 10^{-6}$  &   $1.482\cdot 10^{15}$    &  0.4500\\
		8  &   15 &   $1.668\cdot 10^{-8}$   &   $1.476\cdot 10^{15}$   &   0.6750\\
		15 &   14  &   $1.668\cdot 10^{-8}$  &   $7.965\cdot 10^{2}$    &   3.417\\
		\hline
	\end{tabular}
\end{table}

\begin{figure}[h!]
	\centering
	\includegraphics[width=170mm]{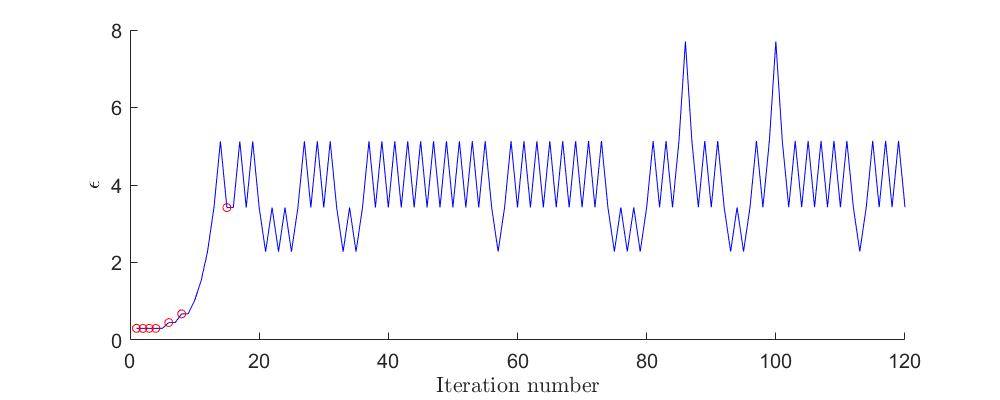}
	\caption{Evolution of $\epsilon$ with the iteration number. Example \ref{ej}.}
	\label{eps}
\end{figure}

As we can see in Figure \ref{eps}, in this example, the sparsification parameter starts with a value of $\epsilon=0.3$, and succeeds at eliminating most edges in the first few iterations. After getting a network with $15$ edges, there are some iterations in which Algorithm \ref{sparsify} does not remove any edge. Therefore, Algorithm \ref{selfhealing} looks for looser sparse approximations of the network increasing the value of $\epsilon$, but decreasing its value when a sparse approximation does not fit the data up to the desired tolerance; until at iteration number $15$. In this iteration the algorithm arrives to the real topology, with almost the real values of the admittance. After this iteration, the algorithm never removes any more edges, and the value of $\epsilon$ tends to oscillate, taking values between $2.28$ and $7.69$. In Figure \ref{eps}, the iterations in which edges were removed are marked with a circle.

}	
\end{example}

\begin{example}\label{hea}
{\rm	In the next experiment, we use a DC network with the topology of the Heawood graph. The results of the application of Algorithm \ref{selfhealing} with an initial value of $\epsilon$ equal to $0.1$ are shown in Table \ref{heawood_tabla}. The algorithm is able to find the topology of the real network in just $5$ iterations in this case.

\begin{figure}[h!]
%\sidecaption
\centering
	\includegraphics[width=95mm]{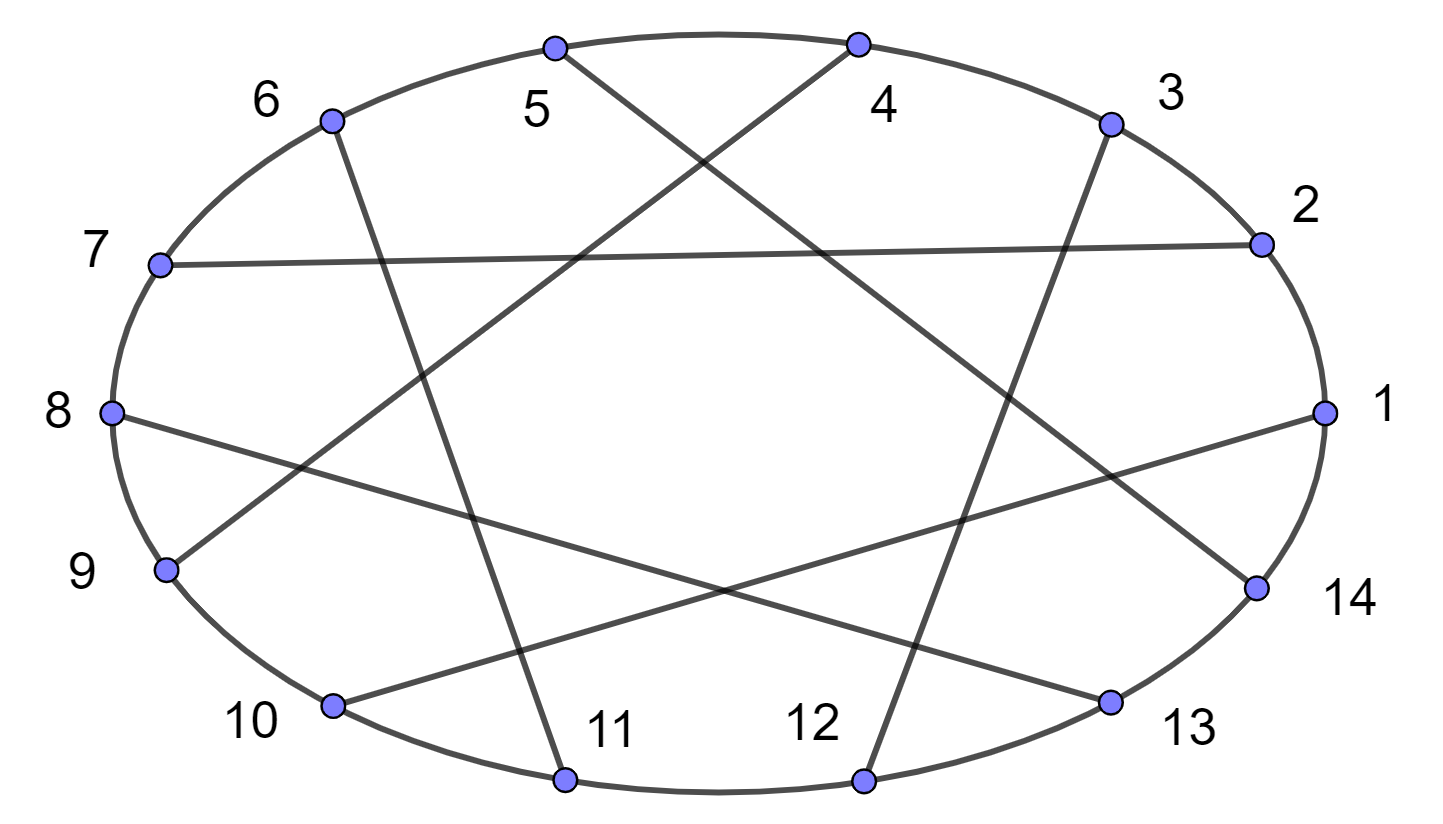}
	\caption{Iteration number: $5$, $|E(\Gamma)|=21$. (Heawood graph network).}
	\label{heawod}
\end{figure}

\begin{table}[h!]
	\centering
	\caption{Heawood graph network fitting results.}
	\label{heawood_tabla}
\begin{tabular}{p{3cm} p{1.5cm} p{2.5cm} p{2.5cm} p{1.5cm}}		\hline
		Iteration number &   $|E(\Gamma)|$  &  $\rms(\Gamma,\bm u,\bm s)$  &  $\kappa(\mathcal{M}_{E(\Gamma),\bm u})$  &  $\epsilon$ \\
		\hline
		1 &   91  &   $3.705\cdot 10^{-6}$   &   $2.119\cdot 10^{4}$   &    0.1000\\
		2  &  35 &   $9.267\cdot 10^{-7}$  &   $2.472\cdot 10^{3}$    &   0.1000\\
		3  &   24  &   $6.378\cdot 10^{-7}$   &   $5.365\cdot 10^{2}$   &   0.1000\\
		4  &   21 &   $9.164\cdot 10^{-10}$  &   $2.586\cdot 10^{2}$    &   0.1000\\
		\hline
	\end{tabular}
\end{table}

In Table \ref{heawood_cota_tabla} we compare the $\rms$ of the sparse approximation in each iteration that removes any edge with the upper bound on this $\rms$ given by Theorem \ref{cota}.

\begin{table}[h!]
	\centering
	\caption{Heawood graph $\rms $ on the sparse approximations.}
	\label{heawood_cota_tabla}
	\begin{tabular}{p{3cm} p{4cm} p{3.5cm}}
		\hline
		Iteration number &Upper bound on the $\rms $ \, \, \, \, \, \,  given by Theorem \ref{cota} & $\rms $ of the sparse \, \, \, \, \, \, \, \,
 approximation\\
		\hline
		2&     $0.4429$   &    0.0369 \\
		3&    $0.4432$  &   0.0308\\
		4&     $0.4432$  &   0.0388\\
		5&      $0.4432$   &   0.0318\\
		\hline
	\end{tabular}
\end{table}

The $\rms$ of the sparse approximations are an order of magnitude smaller than the upper bound guaranteed by Theorem \ref{cota}. The network recovery step that follows any sparsification step decreases the $\rms$ of the network various orders of magnitude. It helps at keeping the $\rms$ of the network below $\tol$ and lets us avoid the possibility of accumulation of error on successive sparsification steps, which could lead to a network that does not correctly fit the data.
}
\end{example}

\begin{example}\label{ejkerber}
{\rm	Finally, we present an example of sparse network recovery on the AC low voltage network called ``Landnetz Freileitung 1 network" (Figure \ref{kerber}), that belongs to the Kerber test networks \cite{ejemplo_ac}.
\begin{figure}[h!]
	\centering
	\includegraphics[width=120mm]{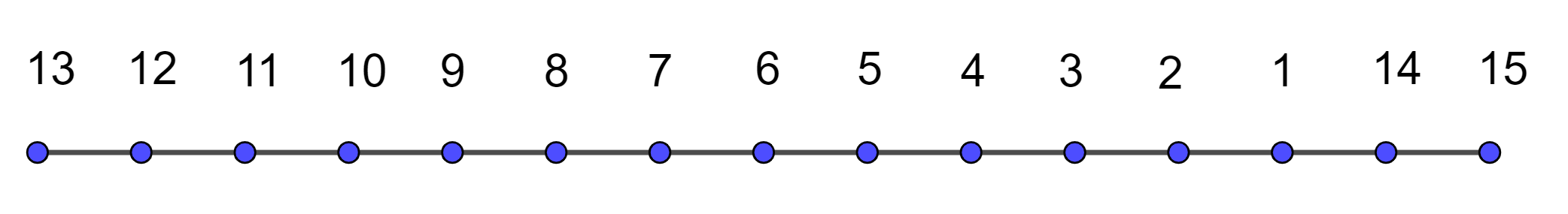}
	\caption{Landnetz Freileitung 1 network.}
	\label{kerber}
\end{figure}

In the data pair  $(\bm u,\bm s)$ there is no injected power at node $x_{14}$ for any component of $\bm s$, that is, $s_j(x_{14})=0$ and thus the potential $u_j$ is harmonic on $x_{14}$, for all $j=1,...,m$. We apply Algorithm \ref{selfhealing} with an initial value of $\epsilon$ equal to $0.3$; which gives the results in Table \ref{kerber_tabla}.

\begin{table}[h!]
	\centering
	\caption{Kerber network fitting results.}
	\label{kerber_tabla}
\begin{tabular}{p{3cm} p{1.5cm} p{2.5cm} p{2.5cm} p{1.5cm}}
		\hline
		Iteration number &   $|E(\Gamma)|$  &  $\rms(\Gamma,\bm u,\bm s)$  &  $\kappa(\mathcal{M}_{E(\Gamma),\bm u})$  &  $\epsilon$ \\
		\hline
		1 &   105 &   $2.069 \cdot 10^{-7}$  &   $1.488 \cdot 10^{15}$ &    0.3 \\
		2  &  16 &   $2.757 \cdot 10^{-7}$ &   $3.829 \cdot 10^{14}$  &   0.3\\
		3  &  15  &   $2.485\cdot 10^{-12}$  &   $3.825 \cdot 10^{14}$ &   0.3\\
		10 &   14  &   $3.885\cdot 10^{-7}$ &   $1.044 \cdot 10^{2}$  &   3.417\\
		34 &   13  &   $2.288\cdot 10^{-12}$ &   $1.029 \cdot 10^{2}$  &   5.126\\
		\hline
	\end{tabular}
\end{table}

The final network $\Gamma'=(V,a')$ has one edge less than the real one $\Gamma=(V,a)$, and has the following topology, which is shown in Figure \ref{kerber_sparse}: the edges $\{x_{1},x_{14}\}$ and $\{x_{14},x_{15}\}$ of the real network are substituted in the recovered network by edge $\{x_{1},x_{15}\}$. The values of admittance on these edges in the real network $\Gamma$ are $a(x_{1}, x_{14})=16.7913-2.6154i$ and $a(x_{14}, x_{15})=1.1999-3.8157i$, and the value of the admittance of $\Gamma'$ at edge $\{x_{1},x_{15}\}$ is $a'(x_{1}, x_{15})=1.6852 - 3.1333i$. The values of $a$ and $a'$ at the edges that are common to both networks, $\Gamma$ and $\Gamma'$, are practically equal.

\begin{figure}[h!]
	\centering
	\includegraphics[width=120mm]{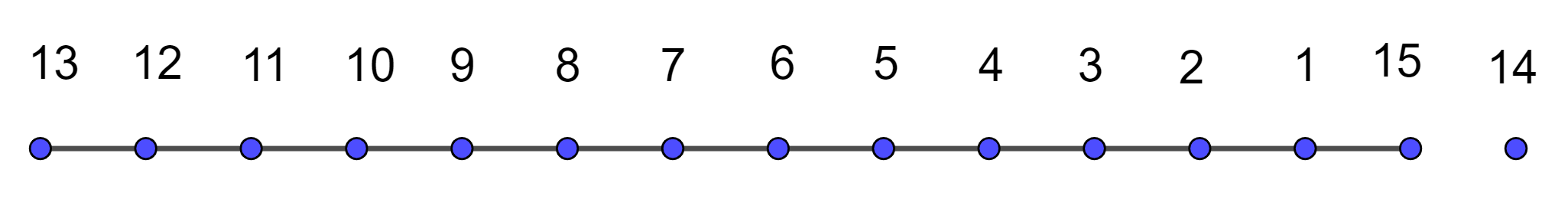}
	\caption{Iteration number: $34$, $|E(\Gamma)|=13$. (Kerber network).}
	\label{kerber_sparse}
\end{figure}

The difference between $\Gamma$ and $\Gamma'$ is analogous to the difference in Example \ref{ej} between the real network of that example whose set of edges is $\left\{e_{xy},e_{yz}\right\}$ and the recovered network of that example whose set of edges is $\left\{e_{xz}\right\}$. That is, it is satisfied that 
\begin{equation}\notag
	\frac{1}{a'(x_{1}, x_{15})}=\frac{1}{a(x_{1}, x_{14})}+\frac{1}{a(x_{14}, x_{15})},
\end{equation}
which is the AC analogue to (\ref{equivalente}). As a consequence, $\Gamma'$ is the union of the isolated vertex $\left\{x_{14}\right\}$ and the Kron reduction of $\Gamma$ with respect to $\left\{x_{14}\right\}$. Also, $a'(x_{1}, x_{15})$ is equal to the effective admittance in $\Gamma$ between $x_{1}$ and $x_{15}$, $a^e(x_{1}, x_{15})$. Analogously to Example \ref{ej}, the Laplacian of $\Gamma'$ is the Dirichlet-to-Neumann map of $\Gamma$ and $\left\{x_{14}\right\}$, which gives us the relationship between potential and injected current at $V\setminus\left\{x_{14}\right\}=\left\{x, z\right\}$ when the potential is harmonic on $x_{14}$.

}
\end{example}

In all the examples, the condition number of the operators
$\mathcal{M}_{E(\Gamma),\bm u}$ decreases dramatically in the algorithm, being much lower in the last electrical network fitted than in the first one. This is especially pronounced in Examples \ref{ej} and \ref{ejkerber}, in which the condition number in the first iteration is of the order of magnitude of $10^{15}$, and it is 13 orders of magnitude lower in the last network topology obtained.

In some cases, the $\rms$ of a network decreases after removing some edges and recovering a network solving Problem \ref{fijo} with the new set of edges. We conjecture that this phenomenon might be due to the finite precision of the interior-point solver.

As a conclusion, the application of Algorithm \ref{selfhealing} to different problems of sparse recovery of an electrical network has yielded promising results. In all cases the algorithm recovers a sparse network after few iterations. In this sparse network topology, the Problem \ref{fijo} of network recovery becomes well-posed.

\section*{Acknowledgments}

The author was funded by a FPI grant of the Research Project PGC2018-
 096446-B-C21 (with the help of the FEDER Program).
 The author was partially supported by the project PID2022-138906NB-C21 funded by
 MICIU/AEI/ 10.13039/501100011033 and by ERDF/EU. The author was partially supported by
 the project 2024PROD00092 funded by AGAUR.

%%%%%%%%%%%%%%%%%%%%%%%%%%%%%%%%%%%%%%%%%%%%%%%%%%%%%%%%%%%%%%%%%%%%%%%%%%%%%%%%%%%%%%%%%%%%%%%%%%%%%%%%%%%%%%%%%%%%%%%%%%%%%%%%%%%%%%%%%%%%%%%%%%%%%%%%%%%%%%%%%%%%%%%%%%%%%%%%%%%%%%%%%%%%%%%%%%%%%%%%%%%%%%%%%%%%%%%%%%%%%%%%%%%%%%%%%%%%%%%%%%%%%%%%%%%%%%%%%%%%%%%%%%%%%%%%%%%%%%%%%%%%%%%%%%%%%%%%%%%%%%%%%%%

%\begin{table}[!t]
%\label{tab:1}  
%\begin{tabular}{p{2cm}p{2.4cm}p{2cm}p{4.9cm}}
%\hline\noalign{\smallskip}
%Classes & Subclass & Length & Action Mechanism  \\
%\noalign{\smallskip}\svhline\noalign{\smallskip}
%Translation & mRNA$^a$  & 22 (19--25) & Translation repression, mRNA cleavage\\
%Translation & mRNA cleavage & 21 & mRNA cleavage\\
%Translation & mRNA  & 21--22 & mRNA cleavage\\
%Translation & mRNA  & 24--26 & Histone and DNA Modification\\
%\end{tabular}
%\end{table}
%

%\input{references}
\bibliographystyle{spmpsci}

\bibliography{biblio}

\end{document}